\documentclass{amsart}

\usepackage{amsmath,amssymb,amsthm,texdraw,color,citesort}
\usepackage{rotating}

\numberwithin{equation}{section}

\newsavebox{\tmpfiga}
\newsavebox{\tmpfigb}
\newsavebox{\tmpfigc}
\newsavebox{\tmpfigd}
\newsavebox{\tmpfige}
\newsavebox{\tmpfigf}
\newsavebox{\tmpfigg}
\newsavebox{\tmpfigh}
\newsavebox{\tmpfigi}
\newsavebox{\tmpfigj}
\newsavebox{\tmpfigk}
\newsavebox{\tmpfigl}
\newsavebox{\tmpfigm}
\newsavebox{\tmpfign}

\theoremstyle{plain}
\newtheorem{thm}{Theorem}[section]
\newtheorem{prop}[thm]{Proposition}
\newtheorem{lem}[thm]{Lemma}

\theoremstyle{definition}
\newtheorem{df}[thm]{Definition}

\newtheorem{example}[thm]{Example}
\theoremstyle{remark}
\newtheorem{remark}[thm]{Remark}

\newcommand{\nc}{\newcommand}

\nc{\eit}{\tilde{e}_i}
\nc{\fit}{\tilde{f}_i}
\nc{\ali}{\alpha_i}
\nc{\csa}{\mathfrak{h}}
\nc{\op}{\oplus}
\nc{\ot}{\otimes}
\nc{\C}{\mathbf{C}}
\nc{\Z}{\mathbf{Z}}
\nc{\half}{\frac{1}{2}}
\nc{\slice}{\mathcal{S}^{(1)}}
\nc{\uq}{U_q}
\nc{\uqp}{U_q'}
\nc{\ftil}{\tilde{f}}
\nc{\etil}{\tilde{e}}
\nc{\newcrystal}{\mathcal{C}^{(1)}}
\nc{\cnone}{C_n^{(1)}}
\nc{\oldcrystal}{\mathcal{B}^{(1)}}
\nc{\defi}[1]{\emph{\textbf{#1}}}
\nc{\veps}{\varepsilon}
\nc{\vphi}{\varphi}
\nc{\La}{\Lambda}
\nc{\la}{\lambda}
\nc{\pathspace}{\mathcal{P}}
\nc{\pwspace}{\mathcal{F}}
\nc{\rpwspace}{\mathcal{Y}}
\nc{\gspath}{\mathbf{b}}
\nc{\path}{\mathbf{p}}
\nc{\spaceiso}{\Psi}
\nc{\spaceosi}{\Phi}
\nc{\gswall}{\mathbf{Y}}
\nc{\wall}{\mathbf{Y}}
\nc{\hwc}{\mathcal{B}}

\DeclareMathOperator{\wt}{wt}
\DeclareMathOperator{\cwt}{\overline{wt}}

\begin{document}

\title[Young walls for $\uq(\cnone)$]
{Young wall realization of\\
Crystal Graphs for $\uq(\cnone)$}

\author{Jin Hong}
\author{Seok-Jin Kang$^*$}
\thanks{$^*$This work was partly supported by KOSEF \# 98-0701-01-5-L
and the Young Scientist Award, Korean Academy of Science and Technology.}
\author{Hyeonmi Lee$^\dagger$}
\thanks{$^\dagger$This work was partly supported by BK21 Project,
Mathematical Sciences Division, Seoul National University.}
\address{Korea Institute for Advanced Study\\
         207-43 Cheongryangri-dong, Dongdaemun-gu\\
         Seoul 130-012, Korea}
\email{jinhong,sjkang,hmlee@kias.re.kr}

\begin{abstract}
We give a realization of crystal graphs for basic representations
of the quantum affine algebra $\uq(\cnone)$ using combinatorics
of Young walls. The notion of splitting blocks plays a crucial
role in the construction of crystal graphs.
\end{abstract}

\maketitle

\section{Introduction}

In~\cite{MR93b:17045} and~\cite{MR90m:17023},
Kashiwara and Lusztig independently
developed the \emph{crystal basis} theory
(or \emph{canonical basis} theory) for
integrable modules over quantum groups associated with
symmetrizable Kac-Moody algebras.
A crystal basis can be viewed as a
basis at $q=0$ and has a structure of colored oriented graph,
called the \emph{crystal graphs,}
defined by \emph{Kashiwara operators.}
The crystal graphs have many nice combinatorial features reflecting
the internal structure of integrable modules over quantum groups.
For example, one can compute the characters of integrable
representations by finding an explicit combinatorial description
of crystal graphs.
Moreover, the crystal graphs have extremely
simple behavior with respect to taking the tensor product.
Thus, the crystal basis theory provides us with a very powerful
combinatorial method of studying the structure of integrable
modules over quantum groups.

Let $\uq(\mathfrak{g})$ be a quantum group associated with a
symmetrizable Kac-Moody algebra $\mathfrak{g}$ and let
$V(\la)$ denote the irreducible highest weight
$\uq(\mathfrak{g})$-module with highest weight $\la$.
One of the most interesting problems in crystal basis theory is to find
explicit realization of the crystal graph $\hwc(\la)$
of $V(\la)$, called the \emph{irreducible highest weight crystal}
with \emph{highest weight} $\la$.
When $\mathfrak{g}$ is a classical
Lie algebra, the crystal graph $\hwc(\la)$ can be
realized as the set of semistandard Young tableaux of a given
shape satisfying certain additional conditions~\cite{MR95c:17025}.

For the quantum affine algebras of type $A_{n-1}^{(1)}$, Misra and
Miwa constructed the crystal graphs $\hwc(\la)$ for
basic representations using colored Young diagrams.
Their idea was
extended to construct crystal graphs for irreducible highest
weight $U_q(A_{n-1}^{(1)})$-modules of arbitrary higher
level~\cite{MR93a:17015}.
The crystal graphs constructed in~\cite{MR93a:17015} and~\cite{MR91j:17021}
can be parameterized by certain \emph{paths}
which arise naturally in the theory
of solvable lattice modules.
Motivated by this observation, Kang,
Kashiwara, Misra, Miwa, Nakashima, and Nakayashiki developed the
theory of \emph{perfect crystals} for quantum affine algebras and gave a
realization of crystal graphs $\hwc(\la)$ over
classical quantum affine algebras of arbitrary higher level in
terms of paths~\cite{MR94j:17013,MR94a:17008,MR93g:17027}.

In~\cite{Kang2000}, Kang introduced the notion of \emph{Young walls}
as a new combinatorial scheme for realizing the crystal graphs for
quantum affine algebras.
The Young walls consist of colored blocks with
various shapes that are built on a given \emph{ground-state wall,} and
they can be viewed as generalizations of colored Young diagrams.
For the classical quantum affine algebras of type $A_n^{(1)}$ ($n\ge 1$),
$A_{2n-1}^{(2)}$ ($n\ge 3$), $D_n^{(1)}$ ($n\ge 4$),
$A_{2n}^{(2)}$ ($n\ge 1$), $D_{n+1}^{(2)}$ ($n\ge 2$),
and $B_n^{(1)}$ ($n\ge 3$), the crystal graphs $\hwc(\la)$ of the basic
representations were realized as the affine crystals consisting of
\emph{reduced proper Young walls.}
However, for the quantum affine
algebras of type $C_n^{(1)}$ ($n\ge 2$), the problem of Young wall
realization of crystal graphs was left open.

The purpose of this paper is to fill in this missing part: we
develop the combinatorics of Young walls for the quantum affine algebras
$\uq(\cnone)$ ($n\ge 2$) and give a realization of crystal graphs
$\hwc(\la)$ for the basic representations as the sets
of reduced proper Young walls.
This case is more difficult to deal
with than the other classical quantum affine algebras, because the
level-$1$ perfect crystals for this case are intrinsically of
level-$2$.
The notion of \emph{splitting blocks} was introduced to
overcome this difficulty.
We believe this notion will play a
crucial role in the construction of higher level irreducible
highest weight crystals for all classical quantum affine algebras.

\section{Quantum group $\uq(\cnone)$ and its level-$1$ perfect crystal}%
\label{sec:02}

We refer the readers to references cited in the introduction,
or to books on quantum groups~\cite{MR1881971,MR96m:17029}
for the basic concepts on quantum groups and crystal bases.
Familiarity with at least the following concepts will be assumed:
quantum group, crystal basis, irreducible highest weight crystal,
(abstract) crystal, perfect crystal, ground state path, $\la$-path,
signature (of a path), path space.

A clear understanding of the Young wall theory for any one of
the affine types developed prior to this work will be immensely helpful
in reading this paper,
although not a logical prerequisite.

Let us fix basic notations here.
\begin{itemize}
\item $\uq = \uq(\cnone)$ : quantum group of type $\cnone$ ($n\geq 2$).
\item $I = \{0,1,\dots, n\}$ : index set.
\item $A = (a_{ij})_{i,j\in I}$ : generalized Cartan matrix of type
      $\cnone$.
\item $P^{\vee} = \left(\op_{i\in I} \Z h_i\right) \op \Z d$ :
      dual weight lattice.
\item $\csa = \C \ot_{\Z} P^{\vee}$ : Cartan subalgebra.
\item $\ali, \delta, \La_i$ : simple root, null root, fundamental weight.
\item $P = \left( \op_{i\in I} \Z\La_i \right) \op \Z\delta$ :
      weight lattice.
\item $e_i, K_i, f_i, q^d$ : generators of $\uq(\cnone)$.
\item $\uqp = \uqp(\cnone)$ : subalgebra of $\uq$
      generated by $e_i, K_i, f_i$ ($i\in I$).
\item $\overline{P} = \op_{i\in I} \Z\La_i$ : classical weight lattice.
\item $\wt, \cwt$ : (affine) weight, classical weight.
\item $\hwc(\La_i)$ : irreducible highest weight crystal of highest
      weight $\La_i$.
\item $\oldcrystal$ : level-$1$ perfect crystal of type $\cnone$.
\item $\eit, \fit$ : Kashiwara operators.
\item $\pathspace(\La_i)$ : set of $\La_i$-paths (with crystal structure).
\end{itemize}

We cite two theorems that are crucial for our work.
The first is the \emph{path realization} of irreducible highest
weight crystals.
\begin{thm}\label{thm:21}
The path space is isomorphic to the irreducible highest weight crystal.
\begin{equation*}
\mathcal{B}(\La_i) \cong \pathspace(\La_i).
\end{equation*}
\end{thm}
A perfect crystal of type $\cnone$, introduced in~\cite{MR95h:17016},
is rewritten in a nice form in~\cite{MR99h:17008}.
We use a special case of this result.
\begin{thm}\label{thm:22}
A level-$1$ perfect crystal of type $\cnone$ is given as follows.
\begin{equation*}
\oldcrystal =
 \Big\{ (x_1,\dots,x_n | \bar{x}_n,\dots, \bar{x}_1) \,\Big\vert\,
x_i, \bar{x}_i \in \Z_{\geq0},\
\textstyle\sum_{i=1}^{n} (x_i + \bar{x}_i) = \textnormal{$0$ or $2$} \Big\}.
\end{equation*}
Let us set $b = (x_1,\dots,x_n | \bar{x}_n,\dots, \bar{x}_1)$.
The action of the Kashiwara operator $\fit$ on $\oldcrystal$ is
given as follows.
\begin{equation*}
\tilde{f}_0 b =
\begin{cases}
(x_1 + 2, x_2, \dots, \bar{x}_2, \bar{x}_1) &
 \textnormal{if $x_1 \geq \bar{x}_1$,}\\
(x_1 + 1, x_2, \dots, \bar{x}_2, \bar{x}_1 - 1) &
 \textnormal{if $x_1 = \bar{x}_1 - 1$,}\\
(x_1, x_2, \dots, \bar{x}_2, \bar{x}_1 - 2) &
 \textnormal{if $x_1 \leq \bar{x}_1 - 2$.}
\end{cases}
\end{equation*}
For $i=1,\dots,n-1$,
\begin{equation*}
\fit b =
\begin{cases}
(x_1, \dots, x_i - 1, x_{i+1} + 1, \dots, \bar{x}_1) &
\textnormal{if $x_{i+1} \geq \bar{x}_{i+1}$,}\\
(x_1, \dots, \bar{x}_{i+1}-1, \bar{x}_{i} + 1, \dots, \bar{x}_1) &
\textnormal{if $x_{i+1} < \bar{x}_{i+1}$.}
\end{cases}
\end{equation*}
And,
\begin{equation*}
\tilde{f}_n b = (x_1, \dots, x_n - 1| \bar{x}_n + 1, \dots, \bar{x}_1).
\end{equation*}
The action of the Kashiwara operator $\eit$ on $\oldcrystal$ is
given as follows.
\begin{equation*}
\tilde{e}_0 b =
\begin{cases}
(x_1 - 2, x_2, \dots, \bar{x}_2, \bar{x}_1) &
 \textnormal{if $x_1 \geq \bar{x}_1 + 2$,}\\
(x_1 - 1, x_2, \dots, \bar{x}_2, \bar{x}_1 + 1) &
 \textnormal{if $x_1 = \bar{x}_1 + 1$,}\\
(x_1, x_2, \dots, \bar{x}_2, \bar{x}_1 + 2) &
 \textnormal{if $x_1 \leq \bar{x}_1$.}
\end{cases}
\end{equation*}
For $i=1,\dots,n-1$,
\begin{equation*}
\eit b =
\begin{cases}
(x_1, \dots, x_i + 1, x_{i+1} - 1, \dots, \bar{x}_1) &
\textnormal{if $x_{i+1} > \bar{x}_{i+1}$,}\\
(x_1, \dots, \bar{x}_{i+1} + 1, \bar{x}_{i} - 1, \dots, \bar{x}_1) &
\textnormal{if $x_{i+1} \leq \bar{x}_{i+1}$.}
\end{cases}
\end{equation*}
And,
\begin{equation*}
\tilde{e}_n b = (x_1, \dots, x_n + 1| \bar{x}_n - 1, \dots, \bar{x}_1).
\end{equation*}
The remaining maps describing a crystal are given below.
\begin{align*}
&\vphi_0 (b) =
  1 - \frac{1}{2}\sum_{i=1}^n (x_i + \bar{x}_i) + (\bar{x}_1 - x_1)_+,\\
&\vphi_i (b) =
  x_i + (\bar{x}_{i+1} - x_{i+1})_+ \quad (i=1,\dots,n),\\
&\vphi_n (b) =
  x_n.\\
&\veps_0 (b) =
  1 - \frac{1}{2}\sum_{i=1}^n (x_i + \bar{x}_i) + (x_1 - \bar{x}_1)_+,\\
&\veps_i (b) =
  \bar{x}_i + ({x}_{i+1} - \bar{x}_{i+1})_+ \quad (i=1,\dots,n),\\
&\veps_n (b) =
  \bar{x}_n.\\
&\cwt(b) =
  \sum_{i=0}^n (\vphi_i(b) - \veps_i(b))\La_i.
\end{align*}
Here, $(x)_+ = \max(0,x)$.
\end{thm}
The perfect crystal for $\uq(C_2^{(1)})$ used in~\cite{MR2001j:17029}
is different from the one given in this theorem.

\begin{example}\label{ex:23}
The following is a drawing of the level-$1$ perfect crystal for
$\uq(C_2^{(1)})$ in the form of the above theorem.
Readers familiar with the crystal bases theory will notice
the $\uq(C_2)$-crystal $\hwc(2\La_1) \subset \hwc(\La_1)\ot\hwc(\La_1)$
in the drawing.
This is what we meant by the \emph{level-$2$ nature} of the perfect
crystal for $\uq(\cnone)$ in the introduction.
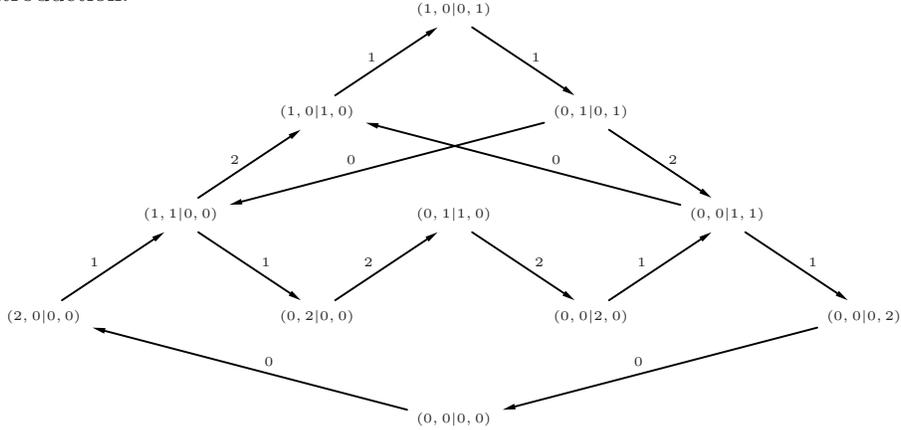
\begin{center}
\begin{texdraw}
\drawdim in \arrowheadsize l:0.065 w:0.03 \arrowheadtype t:F
\fontsize{5}{5}\selectfont
\textref h:C v:C
\drawdim em
\setunitscale 1.9
\htext(0 -6){$(0,0|0,0)$}
\htext(-12 -3){$(2,0|0,0)$}
\htext(-4 -3){$(0,2|0,0)$}
\htext(4 -3){$(0,0|2,0)$}
\htext(12 -3){$(0,0|0,2)$}
\htext(-8 0){$(1,1|0,0)$}
\htext(0 0){$(0,1|1,0)$}
\htext(8 0){$(0,0|1,1)$}
\htext(-4 3){$(1,0|1,0)$}
\htext(4 3){$(0,1|0,1)$}
\htext(0 6){$(1,0|0,1)$}
\move(-11.5 -2.5)\ravec(3 2)
\move(-7.5 0.5)\ravec(3 2)
\move(-3.5 3.5)\ravec(3 2)
\move(-3.5 -2.5)\ravec(3 2)
\move(4.5 -2.5)\ravec(3 2)
\move(4.5 2.5)\ravec(3 -2)
\move(-7.5 -0.5)\ravec(3 -2)
\move(0.5 -0.5)\ravec(3 -2)
\move(0.5 5.5)\ravec(3 -2)
\move(8.5 -0.5)\ravec(3 -2)
\move(6.6 0.3)\ravec(-9.2 2.4)
\move(2.6 2.7)\ravec(-9.2 -2.4)
\move(10.6 -3.3)\ravec(-9.2 -2.4)
\move(-1.4 -5.7)\ravec(-9.2 2.4)
\htext(-10.5 -1.4){$1$}
\htext(-2.5 -1.4){$2$}
\htext(5.5 -1.4){$1$}
\htext(10.5 -1.4){$1$}
\htext(2.5 -1.4){$2$}
\htext(-5.5 -1.4){$1$}
\htext(-6.4 1.6){$2$}
\htext(6.4 1.6){$2$}
\htext(-2.4 4.6){$1$}
\htext(2.4 4.6){$1$}
\htext(-3 1.6){$0$}
\htext(3 1.6){$0$}
\htext(5.4 -4.3){$0$}
\htext(-5.4 -4.3){$0$}
\end{texdraw}
\end{center}
\end{example}

\section{New realization of the level-$1$ perfect crystal}
In this section, we construct the set of \emph{slices}
and obtain a new realization for the level-$1$ perfect crystal $\oldcrystal$.

\subsection{Slices}\label{sec:31}

A \emph{slice} is what will later become a \emph{column}
in our Young walls.

The basic ingredient of our discussion will be
the following colored blocks.
\savebox{\tmpfiga}{\begin{texdraw}
\fontsize{7}{7}\selectfont
\textref h:C v:C
\drawdim em
\setunitscale 1.9
\move(0 0)\lvec(1 0)\lvec(1 0.5)\lvec(0 0.5)\lvec(0 0)
\move(1 0)\rlvec(0.4 0.3)
\move(1 0.5)\rlvec(0.4 0.3)
\move(0 0.5)\rlvec(0.4 0.3)
\move(1.4 0.3)\lvec(1.4 0.8)\lvec(0.4 0.8)
\htext(0.5 0.25){$0$}
\end{texdraw}}%
\savebox{\tmpfigb}{\begin{texdraw}
\fontsize{7}{7}\selectfont
\textref h:C v:C
\drawdim em
\setunitscale 1.9
\move(0 0)\lvec(1 0)\lvec(1 1)\lvec(0 1)\lvec(0 0)
\move(1 0)\rlvec(0.4 0.3)
\move(1 1)\rlvec(0.4 0.3)
\move(0 1)\rlvec(0.4 0.3)
\move(1.4 0.3)\lvec(1.4 1.3)\lvec(0.4 1.3)
\htext(0.5 0.5){$i$}
\end{texdraw}}%
\begin{align*}
\raisebox{-0.25em}{\usebox{\tmpfiga}}
 & \text{\ : half-unit height, unit width, unit depth.}\\
\raisebox{-0.5em}{\usebox{\tmpfigb}}
 & \text{\ : unit height, unit width, unit depth $(i=1,\dots,n)$.}
\end{align*}
To simplify drawings, we shall use just the frontal view
when representing a set of blocks stacked in a wall of unit thickness.
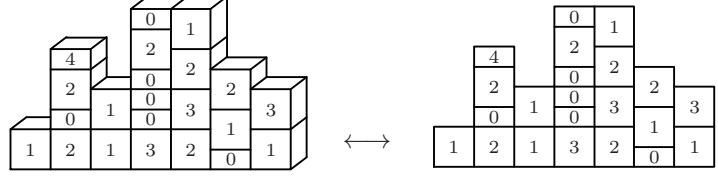
\savebox{\tmpfiga}{\begin{texdraw}
\fontsize{7}{7}\selectfont
\textref h:C v:C
\drawdim em
\setunitscale 1.9
\move(0 0)\lvec(7 0)\lvec(7 2)\lvec(6 2)\lvec(6 2.5)\lvec(5 2.5)
\lvec(5 4)\lvec(3 4)\lvec(3 2)\lvec(2 2)\lvec(2 3)\lvec(1 3)\lvec(1 1)
\lvec(0 1)\lvec(0 0)
\move(1 0)\lvec(1 1)
\move(2 0)\lvec(2 2)
\move(3 0)\lvec(3 2)
\move(4 0)\lvec(4 4)
\move(5 0)\lvec(5 2.5)
\move(6 0)\lvec(6 2)
\move(1 1)\lvec(5 1)
\move(5 0.5)\lvec(6 0.5)
\move(6 1)\lvec(7 1)
\move(1 1.5)\lvec(2 1.5)
\move(1 2.5)\lvec(2 2.5)
\move(3 1.5)\lvec(4 1.5)
\move(3 2)\lvec(4 2)
\move(3 2.5)\lvec(4 2.5)
\move(3 3.5)\lvec(4 3.5)
\move(4 2)\lvec(5 2)
\move(4 3)\lvec(5 3)
\move(5 1.5)\lvec(6 1.5)
\htext(0.5 0.5){$1$}
\htext(1.5 0.5){$2$}
\htext(1.5 1.25){$0$}
\htext(1.5 2){$2$}
\htext(1.5 2.75){$4$}
\htext(2.5 0.5){$1$}
\htext(2.5 1.5){$1$}
\htext(3.5 0.5){$3$}
\htext(3.5 1.25){$0$}
\htext(3.5 1.75){$0$}
\htext(3.5 2.25){$0$}
\htext(3.5 3){$2$}
\htext(3.5 3.75){$0$}
\htext(4.5 0.5){$2$}
\htext(4.5 1.5){$3$}
\htext(4.5 2.5){$2$}
\htext(4.5 3.5){$1$}
\htext(5.5 0.25){$0$}
\htext(5.5 1){$1$}
\htext(5.5 2){$2$}
\htext(6.5 0.5){$1$}
\htext(6.5 1.5){$3$}
\move(7 0)\rlvec(0.4 0.3)
\move(7 1)\rlvec(0.4 0.3)
\move(7 2)\rlvec(0.4 0.3)
\move(6 2)\rlvec(0.4 0.3)
\move(6 2.5)\rlvec(0.4 0.3)
\move(5 2.5)\rlvec(0.4 0.3)
\move(5 3)\rlvec(0.4 0.3)
\move(5 4)\rlvec(0.4 0.3)
\move(4 4)\rlvec(0.4 0.3)
\move(3 4)\rlvec(0.4 0.3)
\move(2 2)\rlvec(0.4 0.3)
\move(2 2.5)\rlvec(0.4 0.3)
\move(2 3)\rlvec(0.4 0.3)
\move(1 3)\rlvec(0.4 0.3)
\move(0 1)\rlvec(0.4 0.3)
\move(7.4 0.3)\lvec(7.4 2.3)\lvec(6.4 2.3)\lvec(6.4 2.8)\lvec(5.4 2.8)
\lvec(5.4 4.3)\lvec(3.4 4.3)
\move(3 2.3)\lvec(2.4 2.3)\lvec(2.4 3.3)\lvec(1.4 3.3)
\move(1 1.3)\lvec(0.4 1.3)
\end{texdraw}}%
\savebox{\tmpfigb}{\begin{texdraw}
\fontsize{7}{7}\selectfont
\textref h:C v:C
\drawdim em
\setunitscale 1.9
\move(0 0)\lvec(7 0)\lvec(7 2)\lvec(6 2)\lvec(6 2.5)\lvec(5 2.5)
\lvec(5 4)\lvec(3 4)\lvec(3 2)\lvec(2 2)\lvec(2 3)\lvec(1 3)\lvec(1 1)
\lvec(0 1)\lvec(0 0)
\move(1 0)\lvec(1 1)
\move(2 0)\lvec(2 2)
\move(3 0)\lvec(3 2)
\move(4 0)\lvec(4 4)
\move(5 0)\lvec(5 2.5)
\move(6 0)\lvec(6 2)
\move(1 1)\lvec(5 1)
\move(5 0.5)\lvec(6 0.5)
\move(6 1)\lvec(7 1)
\move(1 1.5)\lvec(2 1.5)
\move(1 2.5)\lvec(2 2.5)
\move(3 1.5)\lvec(4 1.5)
\move(3 2)\lvec(4 2)
\move(3 2.5)\lvec(4 2.5)
\move(3 3.5)\lvec(4 3.5)
\move(4 2)\lvec(5 2)
\move(4 3)\lvec(5 3)
\move(5 1.5)\lvec(6 1.5)
\htext(0.5 0.5){$1$}
\htext(1.5 0.5){$2$}
\htext(1.5 1.25){$0$}
\htext(1.5 2){$2$}
\htext(1.5 2.75){$4$}
\htext(2.5 0.5){$1$}
\htext(2.5 1.5){$1$}
\htext(3.5 0.5){$3$}
\htext(3.5 1.25){$0$}
\htext(3.5 1.75){$0$}
\htext(3.5 2.25){$0$}
\htext(3.5 3){$2$}
\htext(3.5 3.75){$0$}
\htext(4.5 0.5){$2$}
\htext(4.5 1.5){$3$}
\htext(4.5 2.5){$2$}
\htext(4.5 3.5){$1$}
\htext(5.5 0.25){$0$}
\htext(5.5 1){$1$}
\htext(5.5 2){$2$}
\htext(6.5 0.5){$1$}
\htext(6.5 1.5){$3$}
\end{texdraw}}%
\begin{center}
\raisebox{-0.1em}{\usebox{\tmpfiga}}
\quad\raisebox{0.7em}{$\longleftrightarrow$}\quad\;
\usebox{\tmpfigb}
\end{center}

A set of finitely many blocks, stacked in one column,
following the pattern
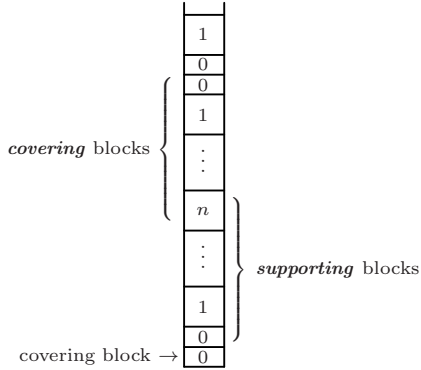
\begin{center}
\begin{texdraw}
\fontsize{7}{7}\selectfont
\textref h:C v:C
\drawdim em
\setunitscale 1.9
\move(0 0)\rlvec(-1 0)
\move(0 0.5)\rlvec(-1 0)
\move(0 1)\rlvec(-1 0)
\move(0 2)\rlvec(-1 0)
\move(0 3.4)\rlvec(-1 0)
\move(0 4.4)\rlvec(-1 0)
\move(0 5.8)\rlvec(-1 0)
\move(0 6.8)\rlvec(-1 0)
\move(0 7.3)\rlvec(-1 0)
\move(0 7.8)\rlvec(-1 0)
\move(0 8.8)\rlvec(-1 0)
\move(0 0)\rlvec(0 9.1)
\move(-1 0)\rlvec(0 9.1)
\htext(-0.5 0.25){$0$}
\htext(-0.5 0.75){$0$}
\htext(-0.5 1.5){$1$}
\vtext(-0.5 2.8){$\cdots$}
\htext(-0.5 3.9){$n$}
\vtext(-0.5 5.2){$\cdots$}
\htext(-0.5 6.3){$1$}
\htext(-0.5 7.05){$0$}
\htext(-0.5 7.55){$0$}
\htext(-0.5 8.3){$1$}
\textref h:L v:C
\htext(0.1 2.5){$\left.\rule{0pt}{3.9em}\right\}$ \defi{supporting} blocks}
\textref h:R v:C
\htext(-1.1 5.5){\defi{covering} blocks $\left\{\rule{0pt}{3.9em}\right.$}
\htext(-1.1 0.25){covering block $\rightarrow$}
\end{texdraw}
\end{center}
is called a \defi{level-$\half$ slice} of type $\cnone$.

We see that, in this repeating pattern,
an $i$-block appears twice in each cycle for $i = 0,\dots,n-1$.
To distinguish the two places, we have given names to these positions or
blocks.
A \defi{covering} block is one that is closer to the $n$-block
that sits below it than to the position for $n$-block above it.
If it is the other way around, it is a \defi{supporting} block.
Notice that, by convention,
each $n$-block is both a supporting block and a covering block.

Any consecutive sequence of blocks in a level-$\half$ slice that contains
one $n$-block and two $i$-blocks for each $i=0,\dots,n-1$
is called a $\delta$.
\begin{remark}
We warn the reader that, even though
$\delta = \alpha_0 + 2\alpha_1 + \cdots + 2\alpha_{n-1} + \alpha_n$
for $\uq(\cnone)$, we are using \emph{two} $0$-blocks
for our definition of $\delta$.
This is because we shall always be using
$0$-blocks in pairs.
For example, in applying $\tilde{f}_0$ action, two $0$-blocks will be added.
\end{remark}

We may \defi{add a $\delta$} to a level-$\half$ slice
or \defi{remove a $\delta$} from a big enough level-$\half$ slice $c$
and write this as $c \pm \delta$.
For example, when dealing with $\uq(C_2^{(1)})$, we have
\begin{center}
\begin{texdraw}
\fontsize{7}{7}\selectfont
\textref h:C v:C
\drawdim em
\setunitscale 1.9
\move(0 0)\lvec(1 0)\lvec(1 2)\lvec(0 2)\lvec(0 0)
\move(0 0.5)\lvec(1 0.5)
\move(0 1)\lvec(1 1)
\htext(0.5 0.25){$0$}
\htext(0.5 0.75){$0$}
\htext(0.5 1.5){$1$}
\end{texdraw}
\raisebox{0.1em}{$\;+\;\delta\ =\ $}
\begin{texdraw}
\fontsize{7}{7}\selectfont
\textref h:C v:C
\drawdim em
\setunitscale 1.9
\move(0 0)\lvec(1 0)\lvec(1 6)\lvec(0 6)\lvec(0 0)
\move(0 0.5)\lvec(1 0.5)
\move(0 1)\lvec(1 1)
\move(0 2)\lvec(1 2)
\move(0 3)\lvec(1 3)
\move(0 4)\lvec(1 4)
\move(0 4.5)\lvec(1 4.5)
\move(0 5)\lvec(1 5)
\htext(0.5 0.25){$0$}
\htext(0.5 0.75){$0$}
\htext(0.5 1.5){$1$}
\htext(0.5 2.5){$2$}
\htext(0.5 3.5){$1$}
\htext(0.5 4.25){$0$}
\htext(0.5 4.75){$0$}
\htext(0.5 5.5){$1$}
\end{texdraw}
\raisebox{0.1em}{.}
\end{center}

\begin{df}
An ordered pair $C = (c_1, c_2)$
of level-$\half$ slices is a
\defi{level-$1$ slice} of type $\cnone$, if
$c_1 \subset c_2 \subset c_1 + \delta$
and if it contains an even number of $0$-blocks.
Each $c_i$ is called the $i$th \defi{layer} of $C$.
The set of all level-$1$ slices is denoted by $\slice$.
\end{df}
We shall often just say \emph{slice}, when dealing with
level-$1$ slices.
Mentally, we picture a level-$1$ slice as two columns with
the first layer placed in front of the second layer, rather than
as an ordered pair.
We explain how to draw a slice with the following example.
\begin{center}
\raisebox{0.3em}{$(\ c_1\;=\ $}
\begin{texdraw}
\fontsize{7}{7}\selectfont
\textref h:C v:C
\drawdim em
\setunitscale 1.9
\move(0 0)\lvec(1 0)\lvec(1 2)\lvec(0 2)\lvec(0 0)
\move(0 0.5)\lvec(1 0.5)
\move(0 1)\lvec(1 1)
\htext(0.5 0.25){$0$}
\htext(0.5 0.75){$0$}
\htext(0.5 1.5){$1$}
\end{texdraw}\;,\quad
\raisebox{0.3em}{$c_2\;=\ $}
\begin{texdraw}
\fontsize{7}{7}\selectfont
\textref h:C v:C
\drawdim em
\setunitscale 1.9
\move(0 0)\lvec(1 0)\lvec(1 4)\lvec(0 4)\lvec(0 0)
\move(0 0.5)\lvec(1 0.5)
\move(0 1)\lvec(1 1)
\move(0 2)\lvec(1 2)
\move(0 3)\lvec(1 3)
\htext(0.5 0.25){$0$}
\htext(0.5 0.75){$0$}
\htext(0.5 1.5){$1$}
\htext(0.5 2.5){$2$}
\htext(0.5 3.5){$1$}
\end{texdraw}
\raisebox{0.3em}{$\ )$}
\qquad\raisebox{0.3em}{$\longleftrightarrow$}\qquad
\raisebox{0.3em}{$C\;=\ $}
\begin{texdraw}
\fontsize{7}{7}\selectfont
\textref h:C v:C
\drawdim em
\setunitscale 1.9
\move(0 0)\lvec(1 0)\lvec(1 2)\lvec(0 2)\lvec(0 0)\ifill f:0.8
\move(0 0)\lvec(1 0)\lvec(1 4)\lvec(0 4)\lvec(0 0)
\move(0 0.5)\lvec(1 0.5)
\move(0 1)\lvec(1 1)
\move(0 2)\lvec(1 2)
\move(0 3)\lvec(1 3)
\htext(0.5 0.25){$0$}
\htext(0.5 0.75){$0$}
\htext(0.5 1.5){$1$}
\htext(0.5 2.5){$2$}
\htext(0.5 3.5){$1$}
\end{texdraw}
\raisebox{0.3em}{$\ =\ $}
\raisebox{-0.4em}{%
\begin{texdraw}
\fontsize{7}{7}\selectfont
\textref h:C v:C
\drawdim em
\setunitscale 1.9 
\move(1 0)\lvec(1 4)\lvec(0 4)\lvec(0 2)\lvec(1 2)
\move(0 3)\lvec(1 3)
\htext(0.5 2.5){$2$}
\htext(0.5 3.5){$1$}
\move(-0.4 -0.3)
\bsegment
\move(0 0)\lvec(1 0)\lvec(1 2)\lvec(0 2)\lvec(0 0)
\move(0 0.5)\lvec(1 0.5)
\move(0 1)\lvec(1 1)
\htext(0.5 0.25){$0$}
\htext(0.5 0.75){$0$}
\htext(0.5 1.5){$1$}
\esegment
\move(0.4 0.3)
\bsegment
\move(1 0)\lvec(1 4)\lvec(0 4)
\esegment
\move(0.6 -0.3)\lvec(1.4 0.3)
\move(0.6 0.2)\lvec(1.4 0.8)
\move(0.6 0.7)\lvec(1.4 1.3)
\move(0.6 1.7)\lvec(1.4 2.3)
\move(-0.4 1.7)\lvec(0 2)
\move(1 3)\lvec(1.4 3.3)
\move(1 4)\lvec(1.4 4.3)
\move(0 4)\lvec(0.4 4.3)
\end{texdraw}}
\end{center}

We now explain the notion of \defi{splitting an $i$-block} in a level-$1$
slice.
Suppose that the top part of some level-$1$ slice $C$
takes one of the following two shapes.
\begin{center}
\raisebox{0.7em}{$C =$}
\begin{texdraw}
\fontsize{7}{7}\selectfont
\textref h:C v:C
\drawdim em
\setunitscale 1.9
\move(0 0)
\bsegment
\move(0 -0.2)\lvec(0 3.1)\lvec(1 3.1)\lvec(1 -0.2)\lvec(0 -0.2)\ifill f:0.8
\esegment
\move(0 4.6)
\bsegment
\move(0 -1.5)\lvec(0 3.3)\lvec(1 3.3)\lvec(1 -1.5)
\move(0 0)\rlvec(1 0)
\move(0 1)\rlvec(1 0)
\move(0 2.3)\rlvec(1 0)
\vtext(0.5 -0.75){$\cdots$}
\htext(0.5 0.5){$n$}
\vtext(0.5 1.65){$\cdots$}
\htext(0.5 2.8){$i$}
\esegment
\move(0 0)
\bsegment
\move(0 -0.2)\lvec(0 3.1)\lvec(1 3.1)\lvec(1 -0.2)
\move(0 0)\rlvec(1 0)
\move(0 0.5)\rlvec(1 0)
\move(0 1)\rlvec(1 0)
\move(0 2.1)\rlvec(1 0)
\htext(0.5 0.25){$0$}
\htext(0.5 0.75){$0$}
\vtext(0.5 1.55){$\cdots$}
\htext(0.5 2.6){$i\!\!-\!\!1$}
\esegment
\end{texdraw}
\qquad\raisebox{0.7em}{or}\qquad
\raisebox{0.7em}{$C =$}
\begin{texdraw}
\fontsize{7}{7}\selectfont
\textref h:C v:C
\drawdim em
\setunitscale 1.9
\move(0 0)
\bsegment
\move(0 -0.2)\lvec(0 3.5)\lvec(1 3.5)\lvec(1 -0.2)\lvec(0 -0.2)\ifill f:0.8
\esegment
\move(0 4.6)
\bsegment
\move(0 -1.1)\lvec(0 3.3)\lvec(1 3.3)\lvec(1 -1.1)
\move(0 0)\rlvec(1 0)
\move(0 0.5)\rlvec(1 0)
\move(0 1)\rlvec(1 0)
\move(0 2.3)\rlvec(1 0)
\vtext(0.5 -0.5){$\cdots$}
\htext(0.5 0.25){$0$}
\htext(0.5 0.75){$0$}
\vtext(0.5 1.65){$\cdots$}
\htext(0.5 2.8){$i\!\!-\!\!1$}
\esegment
\move(0 0)
\bsegment
\move(0 -0.2)\lvec(0 3.5)\lvec(1 3.5)\lvec(1 -0.2)
\move(0 0)\rlvec(1 0)
\move(0 1)\rlvec(1 0)
\move(0 2.6)\rlvec(1 0)
\htext(0.5 0.5){$n$}
\vtext(0.5 1.85){$\cdots$}
\htext(0.5 3){$i$}
\esegment
\end{texdraw}
\end{center}
for some $0 < i < n$.
To \defi{split an $i$-block} in such a level-$1$ slice,
means to break off the top half of the $i$-block and to place it
on top of the $(i-1)$-block, so that it looks like
\begin{center}
\raisebox{0.7em}{$C' =$}
\begin{texdraw}
\fontsize{7}{7}\selectfont
\textref h:C v:C
\drawdim em
\setunitscale 1.9
\move(0 0)
\bsegment
\move(0 -0.2)\lvec(0 3.6)\lvec(1 3.6)\lvec(1 -0.2)\lvec(0 -0.2)\ifill f:0.8
\esegment
\move(0 4.6)
\bsegment
\move(0 -1.5)\lvec(0 2.8)\lvec(1 2.8)\lvec(1 -1.5)
\move(0 -1)\rlvec(1 0)
\move(0 0)\rlvec(1 0)
\move(0 1)\rlvec(1 0)
\move(0 2.3)\rlvec(1 0)
\htext(0.5 -1.25){$i/2$}
\vtext(0.5 -0.45){$\cdots$}
\htext(0.5 0.5){$n$}
\vtext(0.5 1.65){$\cdots$}
\htext(0.5 2.55){$i/2$}
\esegment
\move(0 0)
\bsegment
\move(0 -0.2)\lvec(0 3.1)\lvec(1 3.1)\lvec(1 -0.2)
\move(0 0)\rlvec(1 0)
\move(0 0.5)\rlvec(1 0)
\move(0 1)\rlvec(1 0)
\move(0 2.1)\rlvec(1 0)
\htext(0.5 0.25){$0$}
\htext(0.5 0.75){$0$}
\vtext(0.5 1.55){$\cdots$}
\htext(0.5 2.6){$i\!\!-\!\!1$}
\esegment
\end{texdraw}
\qquad\raisebox{0.7em}{or}\qquad
\raisebox{0.7em}{$C' =$}
\begin{texdraw}
\fontsize{7}{7}\selectfont
\textref h:C v:C
\drawdim em
\setunitscale 1.9
\move(0 0)
\bsegment
\move(0 -0.2)\lvec(0 3)\lvec(1 3)\lvec(1 -0.2)\lvec(0 -0.2)\ifill f:0.8
\esegment
\move(0 4.6)
\bsegment
\move(0 -1.1)\lvec(0 3.8)\lvec(1 3.8)\lvec(1 -1.1)
\move(0 0)\rlvec(1 0)
\move(0 0.5)\rlvec(1 0)
\move(0 1)\rlvec(1 0)
\move(0 2.3)\rlvec(1 0)
\move(0 3.3)\rlvec(1 0)
\vtext(0.5 -0.5){$\cdots$}
\htext(0.5 0.25){$0$}
\htext(0.5 0.75){$0$}
\vtext(0.5 1.65){$\cdots$}
\htext(0.5 2.8){$i\!\!-\!\!1$}
\htext(0.5 3.55){$i/2$}
\esegment
\move(0 0)
\bsegment
\move(0 -0.2)\lvec(0 3.5)\lvec(1 3.5)\lvec(1 -0.2)
\move(0 0)\rlvec(1 0)
\move(0 1)\rlvec(1 0)
\move(0 2.5)\rlvec(1 0)
\move(0 3)\rlvec(1 0)
\htext(0.5 0.5){$n$}
\vtext(0.5 1.85){$\cdots$}
\htext(0.5 2.75){$i/2$}
\esegment
\end{texdraw}
\raisebox{0.7em}{.}
\end{center}
The ``$i/2$'' written in the cut off $i$-blocks are supposed to
convey the idea that this is a half of the $i$-block.

We will never split a $0$-block, but
\defi{splitting an $n$-block} may be done similarly.
\begin{center}
\raisebox{0.7em}{$C =$}
\begin{texdraw}
\fontsize{7}{7}\selectfont
\textref h:C v:C
\drawdim em
\setunitscale 1.9
\move(0 0)
\bsegment
\lvec(0 -0.2)\lvec(0 3.3)\lvec(1 3.3)\lvec(1 -0.2)\lvec(0 -0.2)\ifill f:0.8
\esegment
\move(0 0)
\bsegment
\lvec(0 -0.2)\lvec(0 4.3)\lvec(1 4.3)\lvec(1 -0.2)
\move(0 0)\rlvec(1 0)
\move(0 0.5)\rlvec(1 0)
\move(0 1)\rlvec(1 0)
\move(0 2.3)\rlvec(1 0)
\move(0 3.3)\rlvec(1 0)
\htext(0.5 0.25){$0$}
\htext(0.5 0.75){$0$}
\vtext(0.5 1.65){$\cdots$}
\htext(0.5 2.8){$n\!\!-\!\!1$}
\htext(0.5 3.8){$n$}
\esegment
\end{texdraw}
\qquad\raisebox{0.7em}{$\longmapsto$}\qquad
\raisebox{0.7em}{$C' =$}
\begin{texdraw}
\fontsize{7}{7}\selectfont
\textref h:C v:C
\drawdim em
\setunitscale 1.9
\move(0 0)
\bsegment
\lvec(0 -0.2)\lvec(0 3.8)\lvec(1 3.8)\lvec(1 -0.2)\lvec(0 -0.2)\ifill f:0.8
\esegment
\move(0 0)
\bsegment
\lvec(0 -0.2)\lvec(0 3.8)\lvec(1 3.8)\lvec(1 -0.2)
\move(0 0)\rlvec(1 0)
\move(0 0.5)\rlvec(1 0)
\move(0 1)\rlvec(1 0)
\move(0 2.3)\rlvec(1 0)
\move(0 3.3)\rlvec(1 0)
\htext(0.5 0.25){$0$}
\htext(0.5 0.75){$0$}
\vtext(0.5 1.65){$\cdots$}
\htext(0.5 2.8){$n\!\!-\!\!1$}
\htext(0.5 3.55){$n/2$}
\esegment
\end{texdraw}
\raisebox{0.7em}{.}
\end{center}
\begin{center}
\raisebox{0.7em}{$C =$}
\begin{texdraw}
\fontsize{7}{7}\selectfont
\textref h:C v:C
\drawdim em
\setunitscale 1.9
\move(0 0)
\bsegment
\lvec(0 -0.2)\lvec(0 1)\lvec(1 1)\lvec(1 -0.2)\lvec(0 -0.2)\ifill f:0.8
\esegment
\move(0 0)
\bsegment
\lvec(0 -0.2)\lvec(0 5.6)\lvec(1 5.6)\lvec(1 -0.2)
\move(0 0)\rlvec(1 0)
\move(0 1)\rlvec(1 0)
\move(0 2.3)\rlvec(1 0)
\move(0 2.8)\rlvec(1 0)
\move(0 3.3)\rlvec(1 0)
\move(0 4.6)\rlvec(1 0)
\htext(0.5 0.5){$n$}
\vtext(0.5 1.65){$\cdots$}
\htext(0.5 3.05){$0$}
\htext(0.5 2.55){$0$}
\vtext(0.5 4){$\cdots$}
\htext(0.5 5.1){$n\!\!-\!\!1$}
\esegment
\end{texdraw}
\qquad\raisebox{0.7em}{$\longmapsto$}\qquad
\raisebox{0.7em}{$C' =$}
\begin{texdraw}
\fontsize{7}{7}\selectfont
\textref h:C v:C
\drawdim em
\setunitscale 1.9
\move(0 0)
\bsegment
\lvec(0 -0.2)\lvec(0 0.5)\lvec(1 0.5)\lvec(1 -0.2)\lvec(0 -0.2)\ifill f:0.8
\esegment
\move(0 0)
\bsegment 
\lvec(0 -0.2)\lvec(0 6.1)\lvec(1 6.1)\lvec(1 -0.2)
\move(0 0)\rlvec(1 0)
\move(0 0.5)\rlvec(1 0)
\move(0 1)\rlvec(1 0)
\move(0 2.3)\rlvec(1 0)
\move(0 2.8)\rlvec(1 0)
\move(0 3.3)\rlvec(1 0)
\move(0 4.6)\rlvec(1 0)
\move(0 5.6)\rlvec(1 0)
\htext(0.5 0.25){$n/2$}
\vtext(0.5 1.65){$\cdots$}
\htext(0.5 3.05){$0$}
\htext(0.5 2.55){$0$}
\vtext(0.5 4){$\cdots$}
\htext(0.5 5.1){$n\!\!-\!\!1$}
\htext(0.5 5.85){$n/2$}
\esegment
\end{texdraw}
\raisebox{0.7em}{.}
\end{center}

Simply put, splitting an $i$-block ($i\neq 0$)
is breaking off the top half of
a covering $i$-block and placing it in a supporting $i$-slot.

\begin{remark}
The result obtained after splitting an $i$-block in a slice will
not be considered a level-$1$ slice.
As it will become clearer when we deal with Young walls in
the following sections,
splitting is supposed to be a \emph{temporary} act, used to
see things from a different point of view.
\end{remark}
\begin{remark}\label{rem:34}
If it is possible to split an $i$-block in some slice,
splitting a block of color different from $i$ in the same column is not
possible.
So it makes sense to \emph{split a column if possible.}
\end{remark}

We now explain how to apply some action,
which we denote by $\fit$ ($i=0,\dots,n$),
to the set $\slice \cup \{ 0 \}$.
For $i=1,\dots,n$, we go through the following steps, until we see
a matching case, either to add
one $i$-block to the slice, or to take the result as zero.
\begin{enumerate}
\item
If $i\neq n$ and splitting an $(i+1)$-block is possible,
we take the result to be zero.
\item
If neither of the slots at the top of the two level-$\half$ slices
are for $i$-blocks, the result is zero.
\item
If just one of the two slots is for an $i$-block, place an $i$-block
in the slot.
\item
If we have two $i$-slots at the top, and
either $i=n$ or they are of the same kind (supporting or covering),
do as follows.
\begin{itemize}
\item If they are at different heights, place an $i$-block in
the first layer, i.e., the lower slot.
\item Otherwise, place an $i$-block in the second layer, i.e., the back slot.
\end{itemize}
\item
If we've come this far,
we must have two $i$-slots of different kinds.
Place an $i$-block in the covering slot.
\end{enumerate}
To apply $\ftil_0$, we follow the next steps.
\begin{enumerate}
\item
If it is possible to split a $1$-block, the result is zero.
\item
If there are no $0$-slots available, the result is zero.
\item
If the top of both layers are $0$-slots, place a $0$-block in each
of the two slots.
\item
If only one of the slots is a $0$-slot, place two $0$-blocks in that
slot.
\end{enumerate}

We also define $\eit$ action on $\slice \cup \{0\}$.
For $i\neq 0$, the action of $\eit$ on a slice removes one $i$-block
or sends it to zero following the next set of rules.
\begin{enumerate}
\item
If $i\neq n$ and splitting an $(i+1)$-block is possible,
we take the result to be zero.
\item
If neither of the blocks at the top of the two level-$\half$ slices
are $i$-blocks, the result is zero.
\item
If just one of the two top blocks is an $i$-block, remove the $i$-block.
\item
If the top of both level-$\half$ slices are $i$-blocks,
and if either $i=n$ or they are of the same kind, do as follows.
\begin{itemize}
\item If they are at different heights, remove the $i$-block in
the second layer, i.e., the higher block.
\item Otherwise, remove the $i$-block in the first layer,
i.e., the closer block.
\end{itemize}
\item
We must now have two $i$-blocks of different types at the top of the
two level-$\half$ slices.
Remove the supporting $i$-block.
\end{enumerate}
To apply $\etil_0$, we use the following steps.
\begin{enumerate}
\item
If it is possible to split a $1$-block, the result is zero.
\item
If neither of the two top blocks are $0$-blocks, the result is zero.
\item
If the top of both layers are $0$-blocks, remove a $0$-block from each
of the two layers.
\item
If only one of the top blocks is a $0$-block, remove two $0$-blocks in that
layer.
\end{enumerate}

\subsection{The perfect crystal}

In this subsection, we give a new realization for the level-$1$
perfect crystal of type $\cnone$ by moding out the
repetitive part from the set of slices.

\begin{df}
We may \defi{add a $\delta$} to a slice $C = (c_1, c_2)$,
by changing this into
\begin{equation*}
C + \delta = (c_2, c_1 + \delta).
\end{equation*}
If, for the same slice $C$, we have $c_2 \supset \delta$,
we may also \defi{remove a $\delta$} from $C$,
by changing $C$ into
\begin{equation*}
C - \delta = (c_2 - \delta, c_1).
\end{equation*}
\end{df}

\begin{df}\label{df:35}
Two slices $C$ and $C'$ are \defi{related}, denoted by $C\sim C'$,
if one of the two slices
many be obtained from the other by adding finitely
many $\delta$.
Define
\begin{equation*}
\newcrystal = \slice/\sim.
\end{equation*}
\end{df}

\begin{prop}
The actions $\fit$ and $\eit$,
previously defined on the set of level-$1$ slices,
gives the set $\newcrystal$ a $\uqp(\cnone)$-crystal structure.
\end{prop}
\begin{proof}
Consider the map $C \mapsto C+\delta$, where $C$ denotes a slice.
We may easily check that each of the steps used in defining
$\fit$ and $\eit$ actions on the set of level-$1$ slices commutes
with this map.
So the Kashiwara operators $\fit$ and $\eit$ are well-defined
on $\newcrystal$.

We may now define various other maps as usual.
\begin{align*}
\veps_i ( C ) &= \max\{ n \mid \eit^n C \in \newcrystal \},\\
\vphi_i ( C ) &= \max\{ n \mid \fit^n C \in \newcrystal \},\\
\cwt( C ) &= \sum_i \big(\vphi_i(C) - \veps_i(C)\big) \La_i.
\end{align*}
It is easy to check that these maps define a crystal structure
on $\newcrystal$.
\end{proof}

Recall the finite $\uqp(\cnone)$-crystal $\oldcrystal$ given
in Section~\ref{sec:02}.
We shall define a map from $\oldcrystal$ to $\newcrystal$.
The following preliminary mapping is first needed.
\savebox{\tmpfiga}{\begin{texdraw}
\fontsize{7}{7}\selectfont
\textref h:C v:C
\drawdim em
\setunitscale 1.9
\move(0 -0.2)\lvec(0 0.5)\lvec(1 0.5)\lvec(1 -0.2)
\move(0 0)\lvec(1 0)
\htext(0.5 0.25){$0$}
\end{texdraw}}%
\savebox{\tmpfigb}{\begin{texdraw}
\fontsize{7}{7}\selectfont
\textref h:C v:C
\drawdim em
\setunitscale 1.9
\move(0 -0.2)\lvec(0 1)\lvec(1 1)\lvec(1 -0.2)
\move(0 0)\rlvec(1 0)
\move(0 0.5)\rlvec(1 0)
\htext(0.5 0.25){$0$}
\htext(0.5 0.75){$0$}
\end{texdraw}}%
\savebox{\tmpfigc}{\begin{texdraw}
\fontsize{7}{7}\selectfont
\textref h:C v:C
\drawdim em
\setunitscale 1.9
\move(0 -0.2)\lvec(0 3.5)\lvec(1 3.5)\lvec(1 -0.2)
\move(0 0)\rlvec(1 0)
\move(0 0.5)\rlvec(1 0)
\move(0 1)\rlvec(1 0)
\move(0 2.5)\rlvec(1 0)
\htext(0.5 0.25){$0$}
\htext(0.5 0.75){$0$}
\vtext(0.5 1.75){$\cdots$}
\htext(0.5 3){$i\!\!-\!\!1$}
\end{texdraw}}%
\savebox{\tmpfigd}{\begin{texdraw}
\fontsize{7}{7}\selectfont
\textref h:C v:C
\drawdim em
\setunitscale 1.9
\move(0 -0.2)\lvec(0 3.5)\lvec(1 3.5)\lvec(1 -0.2)
\move(0 0)\rlvec(1 0)
\move(0 0.5)\rlvec(1 0)
\move(0 1)\rlvec(1 0)
\move(0 2.5)\rlvec(1 0)
\htext(0.5 0.25){$0$}
\htext(0.5 0.75){$0$}
\vtext(0.5 1.75){$\cdots$}
\htext(0.5 3){$n$}
\end{texdraw}}%
\savebox{\tmpfige}{\begin{texdraw}
\fontsize{7}{7}\selectfont
\textref h:C v:C
\drawdim em
\setunitscale 1.9
\move(0 -0.2)\lvec(0 6)\lvec(1 6)\lvec(1 -0.2)
\move(0 0)\rlvec(1 0)
\move(0 0.5)\rlvec(1 0)
\move(0 1)\rlvec(1 0)
\move(0 2.5)\rlvec(1 0)
\move(0 3.5)\rlvec(1 0)
\move(0 5)\rlvec(1 0)
\htext(0.5 0.25){$0$}
\htext(0.5 0.75){$0$}
\vtext(0.5 1.75){$\cdots$}
\htext(0.5 3){$n$}
\vtext(0.5 4.25){$\cdots$}
\htext(0.5 5.5){$i$}
\end{texdraw}}%
{\allowdisplaybreaks
\begin{align*}
(0,\dots,0|0,\dots,0) \quad&\longmapsto\quad
 \raisebox{-0.5em}{\usebox{\tmpfiga}}\\[0.3em]
(1,0,\dots,0|0,\dots,0) \quad&\longmapsto\quad
 \raisebox{-0.5em}{\usebox{\tmpfigb}}\\[0.3em]
(0,\dots,0,1,0,\dots,0|0,\dots,0) \quad&\longmapsto\quad
 \raisebox{-0.5em}{\usebox{\tmpfigc}}
 \qquad\text{($1$ at the $i$th place from left)}\\[0.3em]
(0,\dots,0|1,0,\dots,0) \quad&\longmapsto\quad
 \raisebox{-0.8em}{\usebox{\tmpfigd}}\\[0.3em]
(0,\dots,0|0,\dots,0,1,0,\dots,0) \quad&\longmapsto\quad
 \raisebox{-0.8em}{\usebox{\tmpfige}}
 \qquad\text{($1$ at the $i$th place from the right)}
\end{align*}}
Now, to map an element of $\oldcrystal$ to $\newcrystal$,
we first write the element as a sum of two elements,
and map it to $\newcrystal$, using the above preliminary mapping.
The following few examples in the case of $C_2^{(1)}$
should make this clearer.
\savebox{\tmpfiga}{\begin{texdraw}
\fontsize{7}{7}\selectfont
\textref h:C v:C
\drawdim em
\setunitscale 1.9
\move(0 0)\lvec(0 0.5)\lvec(1 0.5)\lvec(1 0)\lvec(0 0)\ifill f:0.8
\move(0 -0.2)\lvec(0 0.5)\lvec(1 0.5)\lvec(1 -0.2)
\move(0 0)\lvec(1 0)
\htext(0.5 0.25){$0$}
\end{texdraw}}%
\savebox{\tmpfigb}{\begin{texdraw}
\fontsize{7}{7}\selectfont
\textref h:C v:C
\drawdim em
\setunitscale 1.9
\move(0 0)\lvec(1 0)\lvec(1 1)\lvec(0 1)\lvec(0 0)\ifill f:0.8
\move(0 -0.2)\lvec(0 3)\lvec(1 3)\lvec(1 -0.2)
\move(0 0)\rlvec(1 0)
\move(0 0.5)\rlvec(1 0)
\move(0 1)\rlvec(1 0)
\move(0 2)\rlvec(1 0)
\htext(0.5 0.25){$0$}
\htext(0.5 0.75){$0$}
\htext(0.5 1.5){$1$}
\htext(0.5 2.5){$2$}
\end{texdraw}}%
\savebox{\tmpfigc}{\begin{texdraw}
\fontsize{7}{7}\selectfont
\textref h:C v:C
\drawdim em
\setunitscale 1.9
\move(0 0)\lvec(1 0)\lvec(1 2)\lvec(0 2)\lvec(0 0)\ifill f:0.8
\move(0 -0.2)\lvec(0 2)\lvec(1 2)\lvec(1 -0.2)
\move(0 0)\rlvec(1 0)
\move(0 0.5)\rlvec(1 0)
\move(0 1)\rlvec(1 0)
\htext(0.5 0.25){$0$}
\htext(0.5 0.75){$0$}
\htext(0.5 1.5){$1$}
\end{texdraw}}%
\begin{alignat*}{2}
(0,0|0,0) &= (0,0|0,0) + (0,0|0,0) &\quad\longmapsto\quad
&\raisebox{-0.5em}{\usebox{\tmpfiga}}\,,\\[0.3em]
(1,0|1,0) &= (1,0|0,0) + (0,0|1,0) &\quad\longmapsto\quad
&\raisebox{-0.5em}{\usebox{\tmpfigb}}\,,\\[0.3em]
(0,2|0,0) &= (0,1|0,0) + (0,1|0,0) &\quad\longmapsto\quad
&\raisebox{-0.5em}{\usebox{\tmpfigc}}\,.\\[0.3em]
\end{alignat*}
Of course, the right hand side should be taken as the equivalence
class in $\slice/\sim$, represented by the drawing.
It is easy to see that this correspondence does not depend on which
of the two summands we decide to map to the first layer.
For example, we could have done
\savebox{\tmpfigb}{\begin{texdraw}
\fontsize{7}{7}\selectfont
\textref h:C v:C
\drawdim em
\setunitscale 1.9
\move(0 0)\lvec(1 0)\lvec(1 3)\lvec(0 3)\lvec(0 0)\ifill f:0.8
\move(0 -0.2)\lvec(0 5)\lvec(1 5)\lvec(1 -0.2)
\move(0 0)\rlvec(1 0)
\move(0 0.5)\rlvec(1 0)
\move(0 1)\rlvec(1 0)
\move(0 2)\rlvec(1 0)
\move(0 3)\rlvec(1 0)
\move(0 4)\rlvec(1 0)
\move(0 4.5)\rlvec(1 0)
\htext(0.5 0.25){$0$}
\htext(0.5 0.75){$0$}
\htext(0.5 1.5){$1$}
\htext(0.5 2.5){$2$}
\htext(0.5 3.5){$1$}
\htext(0.5 4.25){$0$}
\htext(0.5 4.75){$0$}
\end{texdraw}}%
\begin{equation*}
(1,0|1,0) = (0,0|1,0) + (1,0|0,0) \quad\longmapsto\quad
\raisebox{-0.5em}{\usebox{\tmpfigb}}
\end{equation*}
for the second example above.
This might look different at first sight,
but you can check that the drawing on the right hand side belongs
to the same equivalence class as the one given above.

\begin{thm}\label{thm:37}
The map defined above is an isomorphism of $\uqp(\cnone)$-crystals.
\begin{equation*}
\oldcrystal \cong \newcrystal.
\end{equation*}
\end{thm}
\begin{proof}
It is easy to construct an inverse of the above defined map.
Hence the map is bijective.
Checking that the Kashiwara operators of Theorem~\ref{thm:22} and
that of this section is compatible is
also an easy case-by-case comparison.
\end{proof}

\begin{remark}
If we change the definition of a slice slightly to be a column
which extends infinitely downward,
we may state that the set of slices $\slice$ is a realization for the
affinization of $\oldcrystal$.
\end{remark}

We close this section with an example that could be of help in
understanding the proof of Theorem~\ref{thm:37}.
\begin{example}
The following is a drawing of the level-$1$ perfect crystal for
$\uq(C_2^{(1)})$ in the form of $\newcrystal$.
Readers may want to compare this with Example~\ref{ex:23}.
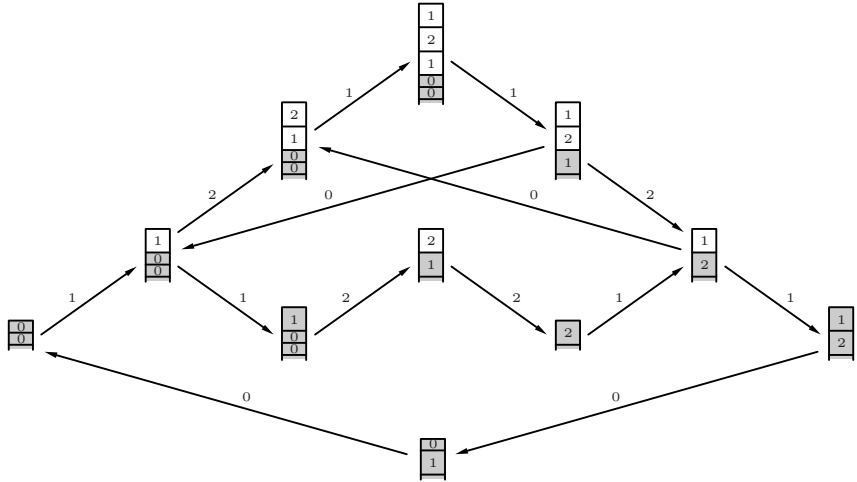
\begin{center}
\begin{texdraw}
\drawdim in \arrowheadsize l:0.065 w:0.03 \arrowheadtype t:F
\fontsize{5}{5}\selectfont
\textref h:C v:C
\drawdim em
\setunitscale 1.9
\move(0 0)
\bsegment
\move(-11.4 -2.5)\ravec(2.8 2)
\move(-7.4 0.5)\ravec(2.8 2)
\move(-3.4 3.5)\ravec(2.8 2)
\move(-3.4 -2.5)\ravec(2.8 2)
\move(4.6 -2.5)\ravec(2.8 2)
\move(4.6 2.5)\ravec(2.8 -2)
\move(-7.4 -0.5)\ravec(2.8 -2)
\move(0.6 -0.5)\ravec(2.8 -2)
\move(0.6 5.5)\ravec(2.8 -2)
\move(8.6 -0.5)\ravec(2.8 -2)
\move(7.3 0)\ravec(-10.6 3)
\move(3.3 3)\ravec(-10.6 -3)
\move(11.3 -3)\ravec(-10.6 -3)
\move(-0.7 -6)\ravec(-10.6 3)
\htext(-10.5 -1.4){$1$}
\htext(-2.5 -1.4){$2$}
\htext(5.5 -1.4){$1$}
\htext(10.5 -1.4){$1$}
\htext(2.5 -1.4){$2$}
\htext(-5.5 -1.4){$1$}
\htext(-6.4 1.6){$2$}
\htext(6.4 1.6){$2$}
\htext(-2.4 4.6){$1$}
\htext(2.4 4.6){$1$}
\htext(-3 1.6){$0$}
\htext(3 1.6){$0$}
\htext(5.4 -4.3){$0$}
\htext(-5.4 -4.3){$0$}
\esegment
\move(-0.35 -0.8)
\bsegment
\move(-8 0)
\bsegment
\setsegscale 0.7
\move(0 -0.15)\lvec(0 1)\lvec(1 1)\lvec(1 -0.15)\lvec(0 -0.15)\ifill f:0.8
\move(0 -0.2)\lvec(0 2)\lvec(1 2)\lvec(1 -0.2)
\move(0 0)\lvec(1 0)
\move(0 0.5)\lvec(1 0.5)
\move(0 1)\lvec(1 1)
\htext(0.5 0.27){$0$}
\htext(0.5 0.77){$0$}
\htext(0.5 1.5){$1$}
\esegment
\move(0 0)
\bsegment
\setsegscale 0.7
\move(0 -0.15)\lvec(0 1)\lvec(1 1)\lvec(1 -0.15)\lvec(0 -0.15)\ifill f:0.8
\move(0 -0.2)\lvec(0 2)\lvec(1 2)\lvec(1 -0.2)
\move(0 0)\lvec(1 0)
\move(0 1)\lvec(1 1)
\htext(0.5 0.5){$1$}
\htext(0.5 1.5){$2$}
\esegment
\move(8 0)
\bsegment
\setsegscale 0.7
\move(0 -0.15)\lvec(0 1)\lvec(1 1)\lvec(1 -0.15)\lvec(0 -0.15)\ifill f:0.8
\move(0 -0.2)\lvec(0 2)\lvec(1 2)\lvec(1 -0.2)
\move(0 0)\lvec(1 0)
\move(0 1)\lvec(1 1)
\htext(0.5 0.5){$2$}
\htext(0.5 1.5){$1$}
\esegment
\esegment
\move(-0.35 -2.8)
\bsegment
\move(-12 0)
\bsegment
\setsegscale 0.7
\move(0 -0.15)\lvec(0 1)\lvec(1 1)\lvec(1 -0.15)\lvec(0 -0.15)\ifill f:0.8
\move(0 -0.2)\lvec(0 1)\lvec(1 1)\lvec(1 -0.2)
\move(0 0)\lvec(1 0)
\move(0 0.5)\lvec(1 0.5)
\htext(0.5 0.27){$0$}
\htext(0.5 0.77){$0$}
\esegment
\move(-4 -0.3)
\bsegment
\setsegscale 0.7
\move(0 -0.15)\lvec(0 2)\lvec(1 2)\lvec(1 -0.15)\lvec(0 -0.15)\ifill f:0.8
\move(0 -0.2)\lvec(0 2)\lvec(1 2)\lvec(1 -0.2)
\move(0 0)\lvec(1 0)
\move(0 0.5)\lvec(1 0.5)
\move(0 1)\lvec(1 1)
\htext(0.5 0.27){$0$}
\htext(0.5 0.77){$0$}
\htext(0.5 1.5){$1$}
\esegment
\move(4 0)
\bsegment
\setsegscale 0.7
\move(0 -0.15)\lvec(0 1)\lvec(1 1)\lvec(1 -0.15)\lvec(0 -0.15)\ifill f:0.8
\move(0 -0.2)\lvec(0 1)\lvec(1 1)\lvec(1 -0.2)
\move(0 0)\lvec(1 0)
\htext(0.5 0.5){$2$}
\esegment
\move(12 -0.3)
\bsegment
\setsegscale 0.7
\move(0 -0.15)\lvec(0 2)\lvec(1 2)\lvec(1 -0.15)\lvec(0 -0.15)\ifill f:0.8
\move(0 -0.2)\lvec(0 2)\lvec(1 2)\lvec(1 -0.2)
\move(0 0)\lvec(1 0)
\move(0 1)\lvec(1 1)
\htext(0.5 0.5){$2$}
\htext(0.5 1.5){$1$}
\esegment
\esegment
\move(-0.3 -6.6)
\bsegment
\setsegscale 0.7
\move(0 -0.15)\lvec(0 1.5)\lvec(1 1.5)\lvec(1 -0.15)\lvec(0 -0.15)\ifill f:0.8
\move(0 -0.2)\lvec(0 1.5)\lvec(1 1.5)\lvec(1 -0.2)
\move(0 0)\lvec(1 0)
\move(0 1)\lvec(1 1)
\htext(0.5 1.27){$0$}
\htext(0.5 0.5){$1$}
\esegment
\move(-0.35 2.2)
\bsegment
\move(-4 0)
\bsegment
\setsegscale 0.7
\move(0 -0.15)\lvec(0 1)\lvec(1 1)\lvec(1 -0.15)\lvec(0 -0.15)\ifill f:0.8
\move(0 -0.2)\lvec(0 3)\lvec(1 3)\lvec(1 -0.2)
\move(0 0)\lvec(1 0)
\move(0 0.5)\lvec(1 0.5)
\move(0 1)\lvec(1 1)
\move(0 2)\lvec(1 2)
\htext(0.5 0.27){$0$}
\htext(0.5 0.77){$0$}
\htext(0.5 1.5){$1$}
\htext(0.5 2.5){$2$}
\esegment
\move(4 0)
\bsegment
\setsegscale 0.7
\move(0 -0.15)\lvec(0 1)\lvec(1 1)\lvec(1 -0.15)\lvec(0 -0.15)\ifill f:0.8
\move(0 -0.2)\lvec(0 3)\lvec(1 3)\lvec(1 -0.2)
\move(0 0)\lvec(1 0)
\move(0 1)\lvec(1 1)
\move(0 2)\lvec(1 2)
\htext(0.5 0.5){$1$}
\htext(0.5 1.5){$2$}
\htext(0.5 2.5){$1$}
\esegment
\esegment
\move(-0.35 4.4)
\bsegment
\setsegscale 0.7
\move(0 -0.15)\lvec(0 1)\lvec(1 1)\lvec(1 -0.15)\lvec(0 -0.15)\ifill f:0.8
\move(0 -0.2)\lvec(0 4)\lvec(1 4)\lvec(1 -0.2)
\move(0 0)\lvec(1 0)
\move(0 0.5)\lvec(1 0.5)
\move(0 1)\lvec(1 1)
\move(0 2)\lvec(1 2)
\move(0 3)\lvec(1 3)
\htext(0.5 0.27){$0$}
\htext(0.5 0.77){$0$}
\htext(0.5 1.5){$1$}
\htext(0.5 2.5){$2$}
\htext(0.5 3.5){$1$}
\esegment
\end{texdraw}
\end{center}
\end{example}
A close study of this example will convince the reader that
choosing to use two half-height blocks (and not a single unit-height block)
for the $\tilde{f}_0$ action was a natural decision.

\section{Young walls}

In this section, we define the set of \emph{reduced proper Young
walls.}
We also define a crystal structure on the set of proper Young walls.

\subsection{Level-$1$ Young walls}

We line up the level-$\half$ slices defined earlier
and consider blocks stacked in the following pattern.
\begin{center}
\begin{texdraw}
\fontsize{7}{7}\selectfont
\textref h:C v:C
\drawdim em
\setunitscale 1.9
\move(0 0)\rlvec(-4.3 0) \move(0 0.5)\rlvec(-4.3 0) \move(0 1)\rlvec(-4.3 0)
\move(0 2)\rlvec(-4.3 0)
\move(0 3.4)\rlvec(-4.3 0)
\move(0 4.4)\rlvec(-4.3 0)
\move(0 5.8)\rlvec(-4.3 0)
\move(0 6.8)\rlvec(-4.3 0)
\move(0 7.3)\rlvec(-4.3 0)
\move(0 7.8)\rlvec(-4.3 0)
\move(0 8.8)\rlvec(-4.3 0)
\move(0 0)\rlvec(0 9.1)
\move(-1 0)\rlvec(0 9.1)
\move(-2 0)\rlvec(0 9.1)
\move(-3 0)\rlvec(0 9.1)
\move(-4 0)\rlvec(0 9.1)
\htext(-0.5 0.25){$0$} \htext(-1.5 0.25){$0$} \htext(-2.5 0.25){$0$}
\htext(-3.5 0.25){$0$} \htext(-0.5 0.75){$0$} \htext(-1.5 0.75){$0$}
\htext(-2.5 0.75){$0$} \htext(-3.5 0.75){$0$}
\htext(-0.5 1.5){$1$} \htext(-1.5 1.5){$1$} \htext(-2.5 1.5){$1$}
\htext(-3.5 1.5){$1$}
\vtext(-0.5 2.7){$\cdots$} \vtext(-1.5 2.7){$\cdots$}
\vtext(-2.5 2.7){$\cdots$} \vtext(-3.5 2.7){$\cdots$}
\htext(-0.5 3.9){$n$} \htext(-1.5 3.9){$n$}
\htext(-2.5 3.9){$n$} \htext(-3.5 3.9){$n$}
\vtext(-0.5 5.1){$\cdots$} \vtext(-1.5 5.1){$\cdots$}
\vtext(-2.5 5.1){$\cdots$} \vtext(-3.5 5.1){$\cdots$}
\htext(-0.5 6.3){$1$} \htext(-1.5 6.3){$1$}
\htext(-2.5 6.3){$1$} \htext(-3.5 6.3){$1$}
\htext(-0.5 7.05){$0$} \htext(-1.5 7.05){$0$} \htext(-2.5 7.05){$0$}
\htext(-3.5 7.05){$0$} \htext(-0.5 7.55){$0$} \htext(-1.5 7.55){$0$}
\htext(-2.5 7.55){$0$} \htext(-3.5 7.55){$0$}
\htext(-0.5 8.3){$1$} \htext(-1.5 8.3){$1$}
\htext(-2.5 8.3){$1$} \htext(-3.5 8.3){$1$}
\end{texdraw}
\end{center}

\begin{df}
A \defi{level-$\half$ weak Young wall} of type $\cnone$ is a set of
blocks, or halves of blocks,
that satisfies the following conditions.
\begin{itemize}
\item It is stacked in the pattern given above.
\item Except for the rightmost column, there is no free
      space to the right of any block.
\end{itemize}
\end{df}

\begin{remark}
In previous works(\cite{Kang2000,MR2001j:17029}),
level-$1$ Young walls were defined to
be built on ground state Young walls.
Implicitly, it was also assumed that the building process was done in
finite steps.
Neither of these conditions are imposed on a level-$\half$ weak Young wall.
For example, the wall
\begin{center}
\begin{texdraw}
\fontsize{7}{7}\selectfont
\textref h:C v:C
\drawdim em
\setunitscale 1.9
\move(-0.2 0)\lvec(7 0)
\move(-0.2 0.5)\lvec(7 0.5)
\move(-0.2 1)\lvec(7 1)
\move(0 0)\lvec(0 1)
\move(1 0)\lvec(1 1)
\move(2 0)\lvec(2 1)
\move(3 0)\lvec(3 2)\lvec(7 2)
\move(4 0)\lvec(4 2)
\move(5 0)\lvec(5 3)\lvec(7 3)
\move(6 0)\lvec(6 3)
\move(7 0)\lvec(7 3)
\htext(0.5 0.25){$0$} \htext(1.5 0.25){$0$} \htext(2.5 0.25){$0$}
\htext(3.5 0.25){$0$} \htext(4.5 0.25){$0$} \htext(5.5 0.25){$0$}
\htext(6.5 0.25){$0$} \htext(0.5 0.75){$0$} \htext(1.5 0.75){$0$}
\htext(2.5 0.75){$0$} \htext(3.5 0.75){$0$} \htext(4.5 0.75){$0$}
\htext(5.5 0.75){$0$} \htext(6.5 0.75){$0$}
\htext(3.5 1.5){$1$} \htext(4.5 1.5){$1$} \htext(5.5 1.5){$1$}
\htext(6.5 1.5){$1$}
\htext(5.5 2.5){$2$} \htext(6.5 2.5){$2$}
\htext(-0.7 1.01){$\cdots$}
\htext(-0.7 0.51){$\cdots$}
\htext(-0.7 0.01){$\cdots$}
\end{texdraw}
\end{center}
is a level-$\half$ weak Young wall of type $C_{n}^{(1)}$,
but not a level-$1$ Young wall in the sense given in previous
works~\cite{MR2001j:17029,Kang2000}.
We also allow the empty wall to be considered a
level-$\half$ weak Young wall.
\end{remark}

For a level-$\half$ weak Young wall $Y$, we define $Y+\delta$ to be the
level-$\half$ weak Young wall obtained by adding a $\delta$ to each and every
column of $Y$.
Here, a $\delta$ is a connected sequence of blocks that contain
one $n$-block and
two $i$-blocks for each $i=0,1,\dots,n-1$.

\begin{df}
An ordered pair $\mathbf{Y} = (Y_1, Y_2)$ is a
\defi{level-$1$ weak Young wall} of type $\cnone$,
if it satisfies the following conditions.
\begin{itemize}
\item Each $Y_i$ is a level-$\half$ weak Young wall.
\item In each column, any halves of blocks for each color add up to form a
      whole block.
      That is, any split blocks come in matching pairs.
\end{itemize}
The level-$\half$ weak Young wall $Y_i$ is called the
\defi{$i$th layer} of the Young wall $\mathbf{Y}$.
\end{df}

\begin{df}
A level-$1$ weak Young wall $\mathbf{Y} = (Y_1, Y_2)$
is a \defi{level-$1$ Young wall}, if it satisfies the following conditions.
\begin{itemize}
\item It contains only whole blocks.
\item $Y_1 \subset Y_2 \subset Y_1 + \delta$.
\item Each column contains an even number of $0$-blocks.
\end{itemize}
\end{df}
In short, a level-$1$ Young wall is a level-$1$ weak Young wall
obtained by concatenating level-$1$ slices.

\subsection{Reduced proper Young walls}

The $i$th column of a level-$1$ Young wall
$\mathbf{Y} = (Y_1, Y_2)$ is denoted by $\mathbf{Y}(i)$.
We choose to number them so that the rightmost column is
named $\mathbf{Y}(0)$.
Note that each column of a level-$1$ Young wall is
a level-$1$ slice.
So we shall utilize the previous notation for drawing level-$1$ slices
when drawing level-$1$ Young walls.
That is, we shall color the first layer gray.
The $i$th column of the $j$th layer $Y_j$ is denoted
by $Y_j(i)$.
We could view the same $Y_j(i)$ also as the $j$th layer of the $i$th column
of $\mathbf{Y}$.

Normally, splitting some block in a column of a level-$1$ Young wall
would not give us a Young wall, or even a weak Young wall.

\begin{df}
Let us be given a level-$1$ Young wall
$\mathbf{Y} = (Y_1, Y_2)$.
The Young wall $\wall$ is \defi{proper} if it satisfies the following
conditions.
\begin{itemize}
\item
When we split every possible column of $\wall$ (see Remark~\ref{rem:34}),
the end result
$\mathbf{Y}' = (Y'_1, Y'_2)$ is a level-$1$ weak Young wall.
\item
For each of the two level-$\half$ weak Young wall $Y'_j$ in the end result,
none of the columns of integer height have the same height.
\end{itemize}
The set of all level-$1$ proper Young walls is denoted by \defi{$\pwspace$}.
\end{df}

Since a column of a Young wall is a slice, we can \emph{add} a $\delta$
to or \emph{remove} a $\delta$ from a column.
To \defi{add a $\delta$} to a column
\begin{equation*}
\mathbf{Y}(i) = (Y_1(i),Y_2(i))
\end{equation*}
means to change this into
\begin{equation*}
\mathbf{Y}(i) + \delta = (Y_2(i), Y_1(i) + \delta).
\end{equation*}
Here, the ``$+\delta$'' on the right hand side should be understood
in the level-$\half$ sense.
Similarly, if $Y_2(i) \supset \delta$,
we may \defi{remove a $\delta$} from the same
column by changing it into
\begin{equation*}
\mathbf{Y}(i) - \delta = (Y_2(i) - \delta, Y_1(i)).
\end{equation*}

\begin{df}
A column in a level-$1$ proper Young wall is said to contain
a \defi{removable $\delta$}, if the Young wall is still proper after
removing a $\delta$ from that column.
A proper Young wall is \defi{reduced}, if none of its columns
contain a removable $\delta$.
The set of all level-$1$ reduced proper Young walls is denoted by
\defi{$\rpwspace$}.
\end{df}

\subsection{The crystal structure}\label{sec:43}

Recall that each column of a Young wall is a slice, so that
we have the actions $\fit$ and $\eit$ defined on
them(Section~\ref{sec:31}).

\begin{df}\hfill
\begin{enumerate}
\item A column in a level-$1$ proper Young wall is $k$ times
   \defi{$i$-admissible}, if $k$ is the maximal number of times
   we may act $\fit$ to the column while remaining
   a proper Young wall.
\item A column in a level-$1$ proper Young wall is $k$ times
   \defi{$i$-removable}, if $k$ is the maximal number of times
   we may act $\eit$ to the column while remaining
   a proper Young wall.
\end{enumerate}
\end{df}

\begin{remark}
Recall that we add two $0$-blocks to a column when applying $\tilde{f}_0$.
So, being $k$ times $0$-admissible will imply that we can place $2k$ number
of $0$-blocks and still obtain a proper Young wall.
\end{remark}
\begin{remark}
Even for $i\neq 0$,
the number of slots in a column of a Young wall,
in which a single $i$-block may be placed while remaining a proper Young wall,
does not necessarily equal the number of times a column is $i$-admissible.
\begin{center}
\begin{texdraw}
\fontsize{7}{7}\selectfont
\textref h:C v:C
\drawdim em
\setunitscale 1.9
\move(-1 -0.15)\lvec(-1 1)\lvec(0 1)\lvec(0 3)\lvec(1 3)\lvec(1 -0.15)
\lvec(-1 -0.15)\ifill f:0.8
\move(-1 -0.2)\lvec(-1 3)\lvec(1 3)
\move(0 -0.2)\lvec(0 4)\lvec(1 4)\lvec(1 -0.2)
\move(-1 2)\lvec(1 2)
\move(-1 1)\lvec(1 1)
\move(-1 0.5)\lvec(1 0.5)
\move(-1 0)\lvec(1 0)
\htext(-0.5 0.25){$0$}
\htext(-0.5 0.75){$0$}
\htext(0.5 0.25){$0$}
\htext(0.5 0.75){$0$}
\htext(-0.5 1.5){$1$}
\htext(-0.5 2.5){$2$}
\htext(0.5 1.5){$1$}
\htext(0.5 2.5){$2$}
\htext(0.5 3.5){$1$}
\htext(-1.6 0.5){$\cdots$}
\htext(1.6 0.5){$\cdots$}
\end{texdraw}
\end{center}
In the above Young wall of type $C_2^{(1)}$,
the left column is only once $1$-admissible.
But we may place a $1$-block in either of the two $i$-slots
in the left column and still obtain a proper Young wall.
\end{remark}

\begin{remark}
The property of a Young wall column being $i$-admissible depends on the
column which sits to the right of the column in consideration.
Likewise, being $i$-removable depends on the left column.
\end{remark}

The action of the Kashiwara operators $\fit$ and $\eit$
on a level-$1$ proper Young wall is defined as follows.
\begin{enumerate}
\item For each column of the Young wall, write under them
      $x$-many 1 followed by $y$-many 0, if the column is $x$ times
      $i$-removable and $y$ times $i$-admissible.
\item From the (half-)infinite list of 0 and 1,
      successively cancel out each $(0,1)$ pair
      to obtain a finite sequence of 1 followed by some 0 (reading
      from left to right).
\item For $\fit$, act $\fit$ to the column corresponding
      to the left-most 0 remaining.
      Set it to zero if no 0 remains, or if
      the action on the column is zero.
\item For $\eit$, act $\eit$ to the column  corresponding
      to the right-most 1 remaining.
      Set it to zero if no 1 remains, or if
      the action on the column is zero.
\end{enumerate}
The $0$ and $1$ placed under the Young wall in the above
process are called \emph{$i$-signature}
of the respective columns or of the Young wall.

It is clear from the definition
that the result obtained after the action of $\fit$ or $\eit$
is still a proper Young wall.
We may now define various other maps as before.
For a proper Young wall $Y$, we define
\begin{align}
\veps_i ( Y ) &= \max\{ n \mid \text{$\eit^n Y$ is nonzero} \},\\
\vphi_i ( Y ) &= \max\{ n \mid \text{$\fit^n Y$ is nonzero} \},\\
\cwt( Y ) &= \sum_i \big(\vphi_i(Y) - \veps_i(Y)\big) \La_i.\label{eq:43}
\end{align}
It is easy to check that these maps give a $\uqp(\cnone)$-crystal
structure to the set of Young walls.

\begin{prop}
The set $\pwspace$ of all level-$1$ proper Young walls
forms an \textup{(}abstract\textup{)}
$\uqp(\cnone)$-crystal.
\end{prop}

\section{Irreducible highest weight crystals}

In this section, we show that the set of all reduced proper Young walls
built on some \emph{ground state wall}
is isomorphic to the irreducible highest weight crystal
of appropriate highest weight.

\subsection{Ground state Young walls}

Let us denote by $\pathspace$, the set of all half-infinite tensor
product of elements from the perfect crystal $\oldcrystal$.
In particular, we have $\pathspace(\La_k) \subset \pathspace$.
By using the composition of maps
\begin{equation*}
\slice \longrightarrow
\newcrystal \stackrel{_\sim}{\longrightarrow}
\oldcrystal
\end{equation*}
on each column of a proper Young wall, we may
define a map
\begin{equation}\label{eq:51}
\spaceiso : \pwspace \longrightarrow \pathspace.
\end{equation}
Now, fix any ground state path $\gspath_{\La_k} \in \pathspace$
and consider its inverse image $\spaceiso^{-1}(\gspath_{\La_k})$.
For example, when dealing with $\uq(C_2^{(1)})$, all of the following
level-$1$ proper Young walls are sent to $\gspath_{\La_0}$ under
the map $\spaceiso$.
\savebox{\tmpfiga}{\begin{texdraw}
\fontsize{7}{7}\selectfont
\textref h:C v:C
\drawdim em
\setunitscale 1.9
\move(-0.2 0)\lvec(5 0)\lvec(5 0.5)\lvec(-0.2 0.5)\lvec(-0.2 0)\ifill f:0.8
\move(-0.25 0)\lvec(5 0)
\move(-0.25 0.5)\lvec(5 0.5)
\move(0 0)\rlvec(0 0.5)
\move(1 0)\rlvec(0 0.5)
\move(2 0)\rlvec(0 0.5)
\move(3 0)\rlvec(0 0.5)
\move(4 0)\rlvec(0 0.5)
\move(5 0)\rlvec(0 0.5)
\htext(0.5 0.25){$0$}
\htext(1.5 0.25){$0$}
\htext(2.5 0.25){$0$}
\htext(3.5 0.25){$0$}
\htext(4.5 0.25){$0$}
\end{texdraw}}%
\savebox{\tmpfigb}{\begin{texdraw}
\fontsize{7}{7}\selectfont
\textref h:C v:C
\drawdim em
\setunitscale 1.9
\move(-0.2 0)\lvec(5 0)\lvec(5 0.5)\lvec(-0.2 0.5)\lvec(-0.2 0)\ifill f:0.8
\move(-0.25 0)\lvec(5 0)
\move(-0.25 0.5)\lvec(5 0.5)
\move(4 1)\rlvec(1 0)
\move(4 2)\rlvec(1 0)
\move(4 3)\rlvec(1 0)
\move(4 4)\rlvec(1 0)
\move(4 4.5)\rlvec(1 0)
\move(0 0)\rlvec(0 0.5)
\move(1 0)\rlvec(0 0.5)
\move(2 0)\rlvec(0 0.5)
\move(3 0)\rlvec(0 0.5)
\move(4 0)\rlvec(0 4.5)
\move(5 0)\rlvec(0 4.5)
\htext(0.5 0.25){$0$}
\htext(1.5 0.25){$0$} 
\htext(2.5 0.25){$0$} 
\htext(3.5 0.25){$0$} 
\htext(4.5 0.25){$0$} 
\htext(4.5 0.75){$0$} 
\htext(4.5 1.5){$1$} 
\htext(4.5 2.5){$2$} 
\htext(4.5 3.5){$1$} 
\htext(4.5 4.25){$0$} 
\end{texdraw}}%
\savebox{\tmpfigc}{\begin{texdraw}
\fontsize{7}{7}\selectfont
\textref h:C v:C
\drawdim em
\setunitscale 1.9
\move(-0.2 0)\lvec(5 0)\lvec(5 0.5)\lvec(-0.2 0.5)\lvec(-0.2 0)\ifill f:0.8
\move(-0.25 0)\lvec(5 0)
\move(-0.25 0.5)\lvec(5 0.5)
\move(3 1)\rlvec(2 0)
\move(3 2)\rlvec(2 0)
\move(3 3)\rlvec(2 0)
\move(3 4)\rlvec(2 0)
\move(3 4.5)\rlvec(2 0)
\move(0 0)\rlvec(0 0.5)
\move(1 0)\rlvec(0 0.5)
\move(2 0)\rlvec(0 0.5)
\move(3 0)\rlvec(0 4.5)
\move(4 0)\rlvec(0 4.5)
\move(5 0)\rlvec(0 4.5)
\htext(0.5 0.25){$0$}
\htext(1.5 0.25){$0$}
\htext(2.5 0.25){$0$}
\htext(3.5 0.25){$0$}
\htext(4.5 0.25){$0$}
\htext(4.5 0.75){$0$}
\htext(4.5 1.5){$1$} 
\htext(4.5 2.5){$2$}
\htext(4.5 3.5){$1$}
\htext(4.5 4.25){$0$}
\htext(3.5 0.75){$0$}
\htext(3.5 1.5){$1$} 
\htext(3.5 2.5){$2$}
\htext(3.5 3.5){$1$}
\htext(3.5 4.25){$0$}
\end{texdraw}}%
\savebox{\tmpfigd}{\begin{texdraw}
\fontsize{7}{7}\selectfont
\textref h:C v:C
\drawdim em
\setunitscale 1.9
\move(-0.2 0)\lvec(5 0)\lvec(5 4.5)\lvec(4 4.5)\lvec(4 0.5)
\lvec(-0.2 0.5)\lvec(-0.2 0)\ifill f:0.8
\move(-0.25 0)\lvec(5 0)
\move(-0.25 0.5)\lvec(5 0.5)
\move(4 1)\rlvec(1 0)
\move(4 2)\rlvec(1 0)
\move(4 3)\rlvec(1 0)
\move(4 4)\rlvec(1 0)
\move(4 4.5)\rlvec(1 0)
\move(0 0)\rlvec(0 0.5)
\move(1 0)\rlvec(0 0.5)
\move(2 0)\rlvec(0 0.5)
\move(3 0)\rlvec(0 0.5)
\move(4 0)\rlvec(0 4.5)
\move(5 0)\rlvec(0 4.5)
\htext(0.5 0.25){$0$}
\htext(1.5 0.25){$0$} 
\htext(2.5 0.25){$0$} 
\htext(3.5 0.25){$0$} 
\htext(4.5 0.25){$0$} 
\htext(4.5 0.75){$0$} 
\htext(4.5 1.5){$1$} 
\htext(4.5 2.5){$2$} 
\htext(4.5 3.5){$1$} 
\htext(4.5 4.25){$0$} 
\end{texdraw}}%
\savebox{\tmpfige}{\begin{texdraw}
\fontsize{7}{7}\selectfont
\textref h:C v:C
\drawdim em
\setunitscale 1.9
\move(-0.2 0)\lvec(5 0)\lvec(5 4.5)\lvec(4 4.5)\lvec(4 0.5)
\lvec(-0.2 0.5)\lvec(-0.2 0)\ifill f:0.8
\move(-0.25 0)\lvec(5 0)
\move(-0.25 0.5)\lvec(5 0.5)
\move(3 1)\rlvec(2 0)
\move(3 2)\rlvec(2 0)
\move(3 3)\rlvec(2 0)
\move(3 4)\rlvec(2 0)
\move(3 4.5)\rlvec(2 0)
\move(0 0)\rlvec(0 0.5)
\move(1 0)\rlvec(0 0.5)
\move(2 0)\rlvec(0 0.5)
\move(3 0)\rlvec(0 4.5)
\move(4 0)\rlvec(0 4.5)
\move(5 0)\rlvec(0 4.5)
\htext(0.5 0.25){$0$}
\htext(1.5 0.25){$0$}
\htext(2.5 0.25){$0$}
\htext(3.5 0.25){$0$}
\htext(4.5 0.25){$0$}
\htext(4.5 0.75){$0$}
\htext(4.5 1.5){$1$}
\htext(4.5 2.5){$2$}
\htext(4.5 3.5){$1$}
\htext(4.5 4.25){$0$}
\htext(3.5 0.75){$0$}
\htext(3.5 1.5){$1$}
\htext(3.5 2.5){$2$}
\htext(3.5 3.5){$1$}
\htext(3.5 4.25){$0$}
\end{texdraw}}%
\savebox{\tmpfigf}{\begin{texdraw}
\fontsize{7}{7}\selectfont
\textref h:C v:C
\drawdim em
\setunitscale 1.9
\move(-0.2 0)\lvec(5 0)\lvec(5 4.5)\lvec(4 4.5)\lvec(4 0.5)
\lvec(-0.2 0.5)\lvec(-0.2 0)\ifill f:0.8
\move(-0.25 0)\lvec(5 0)
\move(-0.25 0.5)\lvec(5 0.5)
\move(2 1)\rlvec(3 0)
\move(2 2)\rlvec(3 0)
\move(2 3)\rlvec(3 0)
\move(2 4)\rlvec(3 0)
\move(2 4.5)\rlvec(3 0)
\move(0 0)\rlvec(0 0.5)
\move(1 0)\rlvec(0 0.5)
\move(2 0)\rlvec(0 4.5)
\move(3 0)\rlvec(0 4.5)
\move(4 0)\rlvec(0 4.5)
\move(5 0)\rlvec(0 4.5)
\htext(0.5 0.25){$0$}
\htext(1.5 0.25){$0$}
\htext(2.5 0.25){$0$}
\htext(3.5 0.25){$0$}
\htext(4.5 0.25){$0$}
\htext(4.5 0.75){$0$}
\htext(4.5 1.5){$1$}
\htext(4.5 2.5){$2$}
\htext(4.5 3.5){$1$}
\htext(4.5 4.25){$0$}
\htext(3.5 0.75){$0$}
\htext(3.5 1.5){$1$}
\htext(3.5 2.5){$2$}
\htext(3.5 3.5){$1$}
\htext(3.5 4.25){$0$}
\htext(2.5 0.75){$0$}
\htext(2.5 1.5){$1$}
\htext(2.5 2.5){$2$}
\htext(2.5 3.5){$1$}
\htext(2.5 4.25){$0$}
\end{texdraw}}%
\savebox{\tmpfigg}{\begin{texdraw}
\fontsize{7}{7}\selectfont
\textref h:C v:C
\drawdim em
\setunitscale 1.9
\move(-0.2 0)\lvec(5 0)\lvec(5 0.5)\lvec(-0.2 0.5)\lvec(-0.2 0)\ifill f:0.8
\move(-0.25 0)\lvec(5 0)
\move(-0.25 0.5)\lvec(5 0.5)
\move(-0.25 1)\lvec(5 1)
\move(-0.25 2)\lvec(5 2)
\move(-0.25 3)\lvec(5 3)
\move(-0.25 4)\lvec(5 4)
\move(-0.25 4.5)\lvec(5 4.5)
\move(0 0)\rlvec(0 4.5)
\move(1 0)\rlvec(0 4.5)
\move(2 0)\rlvec(0 4.5)
\move(3 0)\rlvec(0 4.5)
\move(4 0)\rlvec(0 4.5)
\move(5 0)\rlvec(0 4.5)
\htext(0.5 0.25){$0$}
\htext(1.5 0.25){$0$}
\htext(2.5 0.25){$0$}
\htext(3.5 0.25){$0$}
\htext(4.5 0.25){$0$}
\htext(4.5 0.75){$0$}
\htext(4.5 1.5){$1$}
\htext(4.5 2.5){$2$}
\htext(4.5 3.5){$1$}
\htext(4.5 4.25){$0$}
\htext(3.5 0.75){$0$}
\htext(3.5 1.5){$1$}
\htext(3.5 2.5){$2$}
\htext(3.5 3.5){$1$}
\htext(3.5 4.25){$0$}
\htext(2.5 0.75){$0$}
\htext(2.5 1.5){$1$}
\htext(2.5 2.5){$2$}
\htext(2.5 3.5){$1$}
\htext(2.5 4.25){$0$}
\htext(1.5 0.75){$0$}
\htext(1.5 1.5){$1$}
\htext(1.5 2.5){$2$}
\htext(1.5 3.5){$1$}
\htext(1.5 4.25){$0$}
\htext(0.5 0.75){$0$}
\htext(0.5 1.5){$1$}
\htext(0.5 2.5){$2$}
\htext(0.5 3.5){$1$}
\htext(0.5 4.25){$0$}
\end{texdraw}}%
\savebox{\tmpfigh}{\begin{texdraw}
\fontsize{7}{7}\selectfont
\textref h:C v:C
\drawdim em
\setunitscale 1.9
\move(-0.2 0)\lvec(5 0)\lvec(5 4.5)\lvec(3 4.5)\lvec(3 0.5)
\lvec(-0.2 0.5)\lvec(-0.2 0)\ifill f:0.8
\move(-0.25 0)\lvec(5 0)
\move(-0.25 0.5)\lvec(5 0.5)
\move(-0.25 1)\lvec(5 1)
\move(-0.25 2)\lvec(5 2)
\move(-0.25 3)\lvec(5 3)
\move(-0.25 4)\lvec(5 4)
\move(-0.25 4.5)\lvec(5 4.5)
\move(0 0)\rlvec(0 4.5)
\move(1 0)\rlvec(0 4.5)
\move(2 0)\rlvec(0 4.5)
\move(3 0)\rlvec(0 4.5)
\move(4 0)\rlvec(0 4.5)
\move(5 0)\rlvec(0 4.5)
\htext(0.5 0.25){$0$}
\htext(1.5 0.25){$0$}
\htext(2.5 0.25){$0$}
\htext(3.5 0.25){$0$}
\htext(4.5 0.25){$0$}
\htext(4.5 0.75){$0$}
\htext(4.5 1.5){$1$}
\htext(4.5 2.5){$2$}
\htext(4.5 3.5){$1$}
\htext(4.5 4.25){$0$}
\htext(3.5 0.75){$0$}
\htext(3.5 1.5){$1$}
\htext(3.5 2.5){$2$}
\htext(3.5 3.5){$1$}
\htext(3.5 4.25){$0$}
\htext(2.5 0.75){$0$}
\htext(2.5 1.5){$1$}
\htext(2.5 2.5){$2$}
\htext(2.5 3.5){$1$}
\htext(2.5 4.25){$0$}
\htext(1.5 0.75){$0$}
\htext(1.5 1.5){$1$}
\htext(1.5 2.5){$2$}
\htext(1.5 3.5){$1$}
\htext(1.5 4.25){$0$}
\htext(0.5 0.75){$0$}
\htext(0.5 1.5){$1$}
\htext(0.5 2.5){$2$}
\htext(0.5 3.5){$1$}
\htext(0.5 4.25){$0$}
\end{texdraw}}%
\savebox{\tmpfigi}{\begin{texdraw}
\fontsize{7}{7}\selectfont
\textref h:C v:C
\drawdim em
\setunitscale 1.9
\move(-0.2 0)\lvec(5 0)\lvec(5 4.5)\lvec(-0.2 4.5)\lvec(-0.2 0)\ifill f:0.8
\move(-0.25 0)\lvec(5 0)
\move(-0.25 0.5)\lvec(5 0.5)
\move(-0.25 1)\lvec(5 1)
\move(-0.25 2)\lvec(5 2)
\move(-0.25 3)\lvec(5 3)
\move(-0.25 4)\lvec(5 4)
\move(-0.25 4.5)\lvec(5 4.5)
\move(0 0)\rlvec(0 4.5)
\move(1 0)\rlvec(0 4.5)
\move(2 0)\rlvec(0 4.5)
\move(3 0)\rlvec(0 4.5)
\move(4 0)\rlvec(0 4.5)
\move(5 0)\rlvec(0 4.5)
\htext(0.5 0.25){$0$}
\htext(1.5 0.25){$0$}
\htext(2.5 0.25){$0$}
\htext(3.5 0.25){$0$}
\htext(4.5 0.25){$0$}
\htext(4.5 0.75){$0$}
\htext(4.5 1.5){$1$}
\htext(4.5 2.5){$2$}
\htext(4.5 3.5){$1$}
\htext(4.5 4.25){$0$}
\htext(3.5 0.75){$0$}
\htext(3.5 1.5){$1$}
\htext(3.5 2.5){$2$}
\htext(3.5 3.5){$1$}
\htext(3.5 4.25){$0$}
\htext(2.5 0.75){$0$}
\htext(2.5 1.5){$1$}
\htext(2.5 2.5){$2$}
\htext(2.5 3.5){$1$}
\htext(2.5 4.25){$0$}
\htext(1.5 0.75){$0$}
\htext(1.5 1.5){$1$}
\htext(1.5 2.5){$2$}
\htext(1.5 3.5){$1$}
\htext(1.5 4.25){$0$}
\htext(0.5 0.75){$0$}
\htext(0.5 1.5){$1$}
\htext(0.5 2.5){$2$}
\htext(0.5 3.5){$1$}
\htext(0.5 4.25){$0$}
\end{texdraw}}%
\begin{alignat*}{3}
\usebox{\tmpfiga}&\qquad & \usebox{\tmpfigb}&\qquad
 & \usebox{\tmpfigc}&\\[0.2em]
\usebox{\tmpfigd}&\qquad & \usebox{\tmpfige}&\qquad
 & \usebox{\tmpfigf}&\\[0.2em]
\usebox{\tmpfigg}&\qquad & \usebox{\tmpfigh}&\qquad & \usebox{\tmpfigi}&
\end{alignat*}
We can see that three of these are reduced.

Let us denote by $\gswall_{\La_k} \in \pwspace$ the unique element
of $\spaceiso^{-1}(\gspath_{\La_k})$ that is contained in all other
elements of $\spaceiso^{-1}(\gspath_{\La_k})$.
That is, we take it to be the \emph{smallest} element
of $\spaceiso^{-1}(\gspath_{\La_k})$.
It is clear that each $\gswall_{\La_k}$ is a reduced proper Young wall.
Explicitly, we have
\savebox{\tmpfiga}{\begin{texdraw}
\fontsize{7}{7}\selectfont
\textref h:C v:C
\drawdim em
\setunitscale 1.9
\move(-0.2 0)\lvec(5 0)\lvec(5 0.5)\lvec(-0.2 0.5)\lvec(-0.2 0)\ifill f:0.8
\move(-0.25 0)\lvec(5 0)
\move(-0.25 0.5)\lvec(5 0.5)
\move(0 0)\rlvec(0 0.5)
\move(1 0)\rlvec(0 0.5)
\move(2 0)\rlvec(0 0.5)
\move(3 0)\rlvec(0 0.5)
\move(4 0)\rlvec(0 0.5)
\move(5 0)\rlvec(0 0.5)
\htext(0.5 0.25){$0$}
\htext(1.5 0.25){$0$}
\htext(2.5 0.25){$0$}
\htext(3.5 0.25){$0$}
\htext(4.5 0.25){$0$}
\end{texdraw}}%
\savebox{\tmpfigb}{\begin{texdraw}
\fontsize{7}{7}\selectfont
\textref h:C v:C
\drawdim em
\setunitscale 1.9
\move(-0.2 0)\lvec(5 0)\lvec(5 2.7)\lvec(-0.2 2.7)\lvec(-0.2 0)\ifill f:0.8
\move(-0.25 0)\rlvec(5.25 0)
\move(-0.25 0.5)\rlvec(5.25 0)
\move(-0.25 1)\rlvec(5.25 0)
\move(-0.25 2.7)\rlvec(5.25 0)
\move(-0.25 3.9)\rlvec(5.25 0)
\move(-0.25 4.9)\rlvec(5.25 0)
\move(-0.25 6.6)\rlvec(5.25 0)
\move(0 0)\rlvec(0 6.6)
\move(1 0)\rlvec(0 6.6)
\move(2 0)\rlvec(0 6.6)
\move(3 0)\rlvec(0 6.6)
\move(4 0)\rlvec(0 6.6)
\move(5 0)\rlvec(0 6.6)
\htext(0.5 0.25){$0$}
\htext(0.5 0.75){$0$}
\htext(0.5 2.3){$k\!\!-\!\!1$}
\htext(0.5 4.4){$n$}
\htext(0.5 6.2){$k$}
\vtext(0.5 1.6){$\cdots$}
\vtext(0.5 3.3){$\cdots$}
\vtext(0.5 5.5){$\cdots$}
\htext(1.5 0.25){$0$}
\htext(1.5 0.75){$0$}
\htext(1.5 2.3){$k\!\!-\!\!1$}
\htext(1.5 4.4){$n$}
\htext(1.5 6.2){$k$}
\vtext(1.5 1.6){$\cdots$}
\vtext(1.5 3.3){$\cdots$}
\vtext(1.5 5.5){$\cdots$}
\htext(2.5 0.25){$0$}
\htext(2.5 0.75){$0$}
\htext(2.5 2.3){$k\!\!-\!\!1$}
\htext(2.5 4.4){$n$}
\htext(2.5 6.2){$k$}
\vtext(2.5 1.6){$\cdots$}
\vtext(2.5 3.3){$\cdots$}
\vtext(2.5 5.5){$\cdots$}
\htext(3.5 0.25){$0$}
\htext(3.5 0.75){$0$}
\htext(3.5 2.3){$k\!\!-\!\!1$}
\htext(3.5 4.4){$n$}
\htext(3.5 6.2){$k$}
\vtext(3.5 1.6){$\cdots$}
\vtext(3.5 3.3){$\cdots$}
\vtext(3.5 5.5){$\cdots$}
\htext(4.5 0.25){$0$}
\htext(4.5 0.75){$0$}
\htext(4.5 2.3){$k\!\!-\!\!1$}
\htext(4.5 4.4){$n$}
\htext(4.5 6.2){$k$}
\vtext(4.5 1.6){$\cdots$}
\vtext(4.5 3.3){$\cdots$}
\vtext(4.5 5.5){$\cdots$}
\end{texdraw}}%
\savebox{\tmpfigc}{\begin{texdraw}
\fontsize{7}{7}\selectfont
\textref h:C v:C
\drawdim em
\setunitscale 1.9
\move(-0.2 0)\lvec(5 0)\lvec(5 2.7)\lvec(-0.2 2.7)\lvec(-0.2 0)\ifill f:0.8
\move(-0.25 0)\rlvec(5.25 0)
\move(-0.25 0.5)\rlvec(5.25 0)
\move(-0.25 1)\rlvec(5.25 0)
\move(-0.25 2.7)\rlvec(5.25 0)
\move(-0.25 3.7)\rlvec(5.25 0)
\move(0 0)\rlvec(0 3.7)
\move(1 0)\rlvec(0 3.7)
\move(2 0)\rlvec(0 3.7)
\move(3 0)\rlvec(0 3.7)
\move(4 0)\rlvec(0 3.7)
\move(5 0)\rlvec(0 3.7)
\htext(0.5 0.25){$0$}
\htext(0.5 0.75){$0$}
\htext(0.5 2.3){$n\!\!-\!\!1$}
\htext(0.5 3.2){$n$}
\vtext(0.5 1.6){$\cdots$}
\htext(1.5 0.25){$0$}
\htext(1.5 0.75){$0$}
\htext(1.5 2.3){$n\!\!-\!\!1$}
\htext(1.5 3.2){$n$}
\vtext(1.5 1.6){$\cdots$}
\htext(2.5 0.25){$0$}
\htext(2.5 0.75){$0$}
\htext(2.5 2.3){$n\!\!-\!\!1$}
\htext(2.5 3.2){$n$}
\vtext(2.5 1.6){$\cdots$}
\htext(3.5 0.25){$0$}
\htext(3.5 0.75){$0$}
\htext(3.5 2.3){$n\!\!-\!\!1$}
\htext(3.5 3.2){$n$}
\vtext(3.5 1.6){$\cdots$}
\htext(4.5 0.25){$0$}
\htext(4.5 0.75){$0$}
\htext(4.5 2.3){$n\!\!-\!\!1$}
\htext(4.5 3.2){$n$}
\vtext(4.5 1.6){$\cdots$}
\end{texdraw}}%
\begin{align*}
\gswall_{\La_0} &=\, \raisebox{-0.2em}{\usebox{\tmpfiga}}\,,\\[0.2em]
\gswall_{\La_k} &=\, \raisebox{-0.2em}{\usebox{\tmpfigb}}
                     \quad\text{(for $k=1,\dots,n-1$)},\\[0.2em]
\gswall_{\La_n} &=\, \raisebox{-0.2em}{\usebox{\tmpfigc}}\,.\\[0.2em]
\end{align*}
\begin{df}
The level-$1$ reduced proper Young wall $\gswall_{\La_k}$ is called
the \defi{ground state Young wall} of \defi{weight $\La_k$}.
Any level-$1$ proper Young wall that contains only finitely
many more blocks than $\gswall_{\La_k}$ is said to have been
\defi{built on  $\gswall_{\La_k}$}.
\end{df}

\subsection{Irreducible highest weight crystal}

\begin{df}
The set of all level-$1$ proper Young walls built on $\gswall_{\La_k}$
is denoted by \defi{$\pwspace(\La_k)$}.
And the set of all level-$1$ reduced proper Young walls built on
$\gswall_{\La_k}$ is denoted by \defi{$\rpwspace(\La_k)$}.
\end{df}

\begin{prop}
The set $\pwspace(\La_k)$ of level-$1$ proper Young wall built on
$\gswall_{\La_k}$ forms a $\uq(\cnone)$-crystal.
\end{prop}
\begin{proof}
It is clear that it forms a $\uqp(\cnone)$-subcrystal of $\pwspace$.
So it suffices to give an \emph{affine} weight to each Young wall which
is compatible with other maps.

Recall the previous definition of
classical weight of a Young wall given by~\eqref{eq:43}.
We may define the affine weight of a Young wall
$\wall\in \pwspace(\La_k)$ by
\begin{equation}
\wt(\wall) = \cwt(\wall)
 - \half
   \text{(number of $0$-blocks in $\wall \setminus \gswall_{\La_k}$)}\,
   \delta.
\end{equation}
In particular, the affine weight of $\gswall_{\La_k}$ is $\La_k$.
It is easy to check that $\pwspace(\La_k)$ forms a $\uq(\cnone)$-crystal
with this definition of affine weights.
\end{proof}

We claim that the map $\spaceiso$ given by~\eqref{eq:51}
is bijective when we restrict
the domain and range as
\begin{equation}\label{eq:52}
\spaceiso : \rpwspace(\La_k) \longrightarrow \pathspace(\La_k).
\end{equation}
To check the surjectivity, we may explicitly construct a reduced
proper Young wall which maps to the path in question, and in the process,
we will notice that the condition \emph{reduced}
forces it to be chosen uniquely.
Here is our main theorem,
which shows that this bijection is a crystal isomorphism.

\begin{thm}\label{thm:54}
The set $\rpwspace(\La_k)$ is a $\uq(\cnone)$-subcrystal of $\pwspace(\La_k)$.
It is isomorphic to the irreducible highest weight crystal
$\mathcal{B}(\La_k)$.
\end{thm}

The rest of this section is devoted mostly to proving this theorem.
Let us denote by $\spaceosi$ the inverse of the map~\eqref{eq:52},
but with the range enlarged as
\begin{equation*}
\spaceosi : \pathspace(\La_k) \longrightarrow \pwspace(\La_k).
\end{equation*}
To prove the theorem,
it suffices to show that this is a strict crystal morphism.
Then the image, which we've already seen to be $\rpwspace(\La_k)$,
would be a subcrystal of $\pwspace(\La_k)$.
The second statement follows from Theorem~\ref{thm:21}.
We shall focus our efforts on
showing that this map commutes with the Kashiwara operator $\fit$.
Other parts of the proof are similar or easy.

Let us first review the action of the Kashiwara operator $\fit$ on
a path element
\begin{equation*}
\path = \cdots \ot p(j) \ot \cdots \ot p(1)\ot p(0) \in \pathspace(\La_k).
\end{equation*}
\begin{enumerate}
\item Under each $p(j)$, write
      $\veps_i(p(j))$-many 1 followed by $\vphi_i(p(j))$-many 0.
\item From the (half-)infinite list of 0 and 1,
      successively cancel out each $(0,1)$ pair
      to obtain a finite sequence of 1 followed by some 0 (reading
      from left to right).
\item Act $\fit$ to the $p(j)$ corresponding to the left-most 0 remaining.
      Set it to zero if no 0 remains, or if the action on $p(j)$ is zero.
\end{enumerate}
This is quite similar to the action of $\fit$ on Young walls defined
in Section~\ref{sec:43}.
Now, recalling the definition of $\spaceosi$ and Theorem~\ref{thm:37},
we find that, to prove Theorem~\ref{thm:54},
it suffices to prove the following lemma.

\begin{lem}\label{lem:55}
The Kashiwara operator $\fit$ acts on the $j$th tensor component
of a path $\path$ if and only if
it acts on the $j$th column of the Young wall $\spaceosi(\path)$.
\end{lem}

Let us start explaining how to prove this with some examples.
Consider the following part of a path.
\begin{equation*}
\cdots \ot (0,\dots,0|0,\dots,0,2) \ot (2,0,\dots,0|0,\dots,0) \ot \cdots
\end{equation*}
Suppose we are dealing with $i=0$ case.
The signature that should be under the left element $(0,\dots,0,2)$ is $00$
and that for the right element $(2,0,\dots,0)$ is $11$.
After canceling out the $(0,1)$ pairs, we are left with \emph{nothing.}

Now consider columns of the reduced proper Young wall which corresponds to
this path under the map $\spaceosi$.
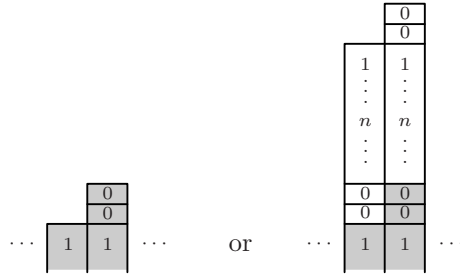
\begin{center}
\begin{texdraw}
\fontsize{7}{7}\selectfont
\textref h:C v:C
\drawdim em
\setunitscale 1.9
\move(1 -0.15)\lvec(1 2)\lvec(0 2)\lvec(0 1)\lvec(-1 1)\lvec(-1 -0.15)
\lvec(1 -0.15)\ifill f:0.8
\move(-1, -0.2)\lvec(-1 1)\lvec(1 1)
\move(1 -0.2)\lvec(1 2)\lvec(0 2)\lvec(0 -0.2)
\move(1 1.5)\lvec(0 1.5)
\htext(0.5 1.25){$0$}
\htext(0.5 1.75){$0$}
\htext(0.5 0.5){$1$}
\htext(-0.5 0.5){$1$}
\htext(1.7 0.5){$\cdots$}
\htext(-1.6 0.5){$\cdots$}
\move(-2 0)\move(2 0)
\end{texdraw}
\qquad\raisebox{0.8em}{or}\qquad
\begin{texdraw}
\fontsize{7}{7}\selectfont
\textref h:C v:C
\drawdim em
\setunitscale 1.9
\move(1 -0.15)\lvec(1 2)\lvec(0 2)\lvec(0 1)\lvec(-1 1)\lvec(-1 -0.15)
\lvec(1 -0.15)\ifill f:0.8
\move(-1 -0.2)\lvec(-1 5.5)\lvec(1 5.5)
\move(0 -0.2)\lvec(0 6.5)\lvec(1 6.5)\lvec(1 -0.2)
\move(-1 1)\lvec(1 1)
\move(-1 1.5)\lvec(1 1.5)
\move(-1 2)\lvec(1 2)
\move(0 6)\lvec(1 6)
\htext(-0.5 0.5){$1$}
\htext(-0.5 1.27){$0$}
\htext(-0.5 1.77){$0$}
\htext(-0.5 3){$\vdots$}
\htext(-0.5 3.55){$n$}
\htext(-0.5 4.5){$\vdots$}
\htext(-0.5 5){$1$}
\htext(0.5 0.5){$1$}
\htext(0.5 1.27){$0$}
\htext(0.5 1.77){$0$}
\htext(0.5 3){$\vdots$}
\htext(0.5 3.55){$n$}
\htext(0.5 4.5){$\vdots$}
\htext(0.5 5){$1$}
\htext(0.5 5.77){$0$}
\htext(0.5 6.27){$0$}
\htext(1.7 0.5){$\cdots$}
\htext(-1.6 0.5){$\cdots$}
\move(-2 0)\move(2 0)
\end{texdraw}
\end{center}
The two drawings are $\delta$-shifts of each other and it does not
matter which of the two drawings we use.
When dealing with the case $i=0$, under the left column we would
write $0$ and under the right column we would write $1$.
After $(0,1)$-pair cancellation, we are again left with \emph{nothing.}
The signatures under path description and the Young wall description
agree after canceling out $(0,1)$ pairs.

We give one more example which is a bit more complicated.
\begin{equation*}
\cdots \ot (0,\dots,0|0,\dots,0) \ot (0,\dots,0|0,\dots,0) \ot \cdots
\end{equation*}
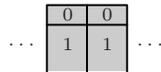
\begin{center}
\begin{texdraw}
\fontsize{7}{7}\selectfont
\textref h:C v:C
\drawdim em
\setunitscale 1.9
\move(1 -0.15)\lvec(1 1.5)\lvec(-1 1.5)\lvec(-1 -0.15)
\lvec(1 -0.15)\ifill f:0.8
\move(-1, -0.2)\lvec(-1 1.5)\lvec(1 1.5)\lvec(1 -0.2)
\move(-1 1)\lvec(1 1)
\move(0 -0.2)\lvec(0 1.5)
\htext(0.5 1.25){$0$}
\htext(-0.5 1.25){$0$}
\htext(0.5 0.5){$1$}
\htext(-0.5 0.5){$1$}
\htext(1.6 0.5){$\cdots$}
\htext(-1.6 0.5){$\cdots$}
\end{texdraw}
\end{center}
As before, the reader is free to use a $\delta$-shift of this Young wall.
When dealing with $\tilde{f}_0$,
the signatures to be written under them are given in the following
table.
\begin{center}
\begin{tabular}{ccccc}
\multicolumn{5}{c}{path}\\
\multicolumn{2}{c}{left} &\quad& \multicolumn{2}{c}{right}\\
$\veps$ & $\vphi$ && $\veps$ & $\vphi$\\
$1$ & $0$ && $1$ & $0$
\end{tabular}
\qquad
\begin{tabular}{ccccc}
\multicolumn{5}{c}{Young wall}\\
\multicolumn{2}{c}{left} &\quad& \multicolumn{2}{c}{right}\\
$\veps$ & $\vphi$ && $\veps$ & $\vphi$\\
? & $\cdot$ && $\cdot$ & ?
\end{tabular}
\end{center}
The first question mark in the above Young wall table signifies that
the number of $1$ that should be written there depends on the column that
sits to left of the left column.
Likewise, the second question mark is to signify that the number of
$0$ to be written depends on the column that comes to its right.
The two dots imply that no $0$ and $1$, respectively,
should be written there.

So in this case, we do not know the complete signature to be written
under the Young wall columns.
Hence, a straightforward comparison of signatures after cancellations
of $(0,1)$ pairs is not possible.
But still, we can verify that what is left of the
left-$\vphi$ and right-$\veps$ signatures,
after the $(0,1)$-pair cancellation,
is the same for the path and Young wall.
In this example, they both amount to nothing.

We can now easily convince ourselves that, to prove Lemma~\ref{lem:55},
it suffices to check exactly this
kind of signature matching between all possible left-right pairs of
perfect crystal elements and their corresponding Young wall columns.
(The right-most column may be dealt with in a similar way.)

We shall do a sketch of this checking now.
Let us deal with the $i=0$ case first.
The following notation will be used to denote various columns of
Young walls.

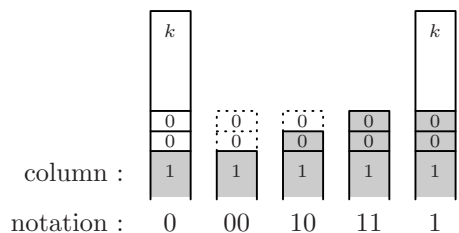
\savebox{\tmpfiga}{\begin{texdraw}
\fontsize{7}{7}\selectfont
\textref h:C v:C
\drawdim em
\setunitscale 1.9
\move(0 -0.15)\lvec(0 1)\lvec(1 1)\lvec(1 -0.15)\lvec(0 -0.15)\ifill f:0.8
\move(0 -0.2)\lvec(0 1)\lvec(1 1)\lvec(1 -0.2)
\lpatt(0.03 0.17)
\move(0 1)\lvec(0 2)\lvec(1 2)\lvec(1 1)
\move(0 1.5)\lvec(1 1.5)
\htext(0.5 0.5){$1$}
\htext(0.5 1.25){$0$}
\htext(0.5 1.75){$0$}
\end{texdraw}}
\savebox{\tmpfigb}{\begin{texdraw}
\fontsize{7}{7}\selectfont
\textref h:C v:C
\drawdim em
\setunitscale 1.9
\move(0 -0.15)\lvec(0 1.5)\lvec(1 1.5)\lvec(1 -0.15)
\lvec(0 -0.15)\ifill f:0.8
\move(0 -0.2)\lvec(0 1.5)\lvec(1 1.5)\lvec(1 -0.2)
\move(0 1)\lvec(1 1)
\lpatt(0.03 0.17)
\move(0 1.5)\lvec(0 2)\lvec(1 2)\lvec(1 1.5)
\htext(0.5 0.5){$1$}
\htext(0.5 1.25){$0$}
\htext(0.5 1.75){$0$}
\end{texdraw}}
\savebox{\tmpfigc}{\begin{texdraw}
\fontsize{7}{7}\selectfont
\textref h:C v:C
\drawdim em
\setunitscale 1.9
\move(0 -0.15)\lvec(0 2)\lvec(1 2)\lvec(1 -0.15)
\lvec(0 -0.15)\ifill f:0.8
\move(0 -0.2)\lvec(0 2)\lvec(1 2)\lvec(1 -0.2)
\move(0 1)\lvec(1 1)
\move(0 1.5)\lvec(1 1.5)
\htext(0.5 0.5){$1$}
\htext(0.5 1.25){$0$}
\htext(0.5 1.75){$0$}
\end{texdraw}}
\savebox{\tmpfigd}{\begin{texdraw}
\fontsize{7}{7}\selectfont
\textref h:C v:C
\drawdim em
\setunitscale 1.9
\move(0 -0.15)\lvec(0 1)\lvec(1 1)\lvec(1 -0.15)\lvec(0 -0.15)\ifill f:0.8
\move(0 -0.2)\lvec(0 4.5)\lvec(1 4.5)\lvec(1 -0.2)
\move(0 1)\lvec(1 1)
\move(0 1.5)\lvec(1 1.5)
\move(0 2)\lvec(1 2)
\htext(0.5 0.5){$1$}
\htext(0.5 1.25){$0$}
\htext(0.5 1.75){$0$}
\htext(0.5 4){$k$}
\end{texdraw}}
\savebox{\tmpfige}{\begin{texdraw}
\fontsize{7}{7}\selectfont
\textref h:C v:C
\drawdim em
\setunitscale 1.9
\move(0 -0.15)\lvec(0 2)\lvec(1 2)\lvec(1 -0.15)\lvec(0 -0.15)\ifill f:0.8
\move(0 -0.2)\lvec(0 4.5)\lvec(1 4.5)\lvec(1 -0.2)
\move(0 1)\lvec(1 1)
\move(0 1.5)\lvec(1 1.5)
\move(0 2)\lvec(1 2)
\htext(0.5 0.5){$1$}
\htext(0.5 1.25){$0$}
\htext(0.5 1.75){$0$}
\htext(0.5 4){$k$}
\end{texdraw}}

\begin{center}
\begin{tabular}{rccccc}
\raisebox{0.7em}{column :} &
\usebox{\tmpfigd}&
\usebox{\tmpfiga}&
\usebox{\tmpfigb}&
\usebox{\tmpfigc}&
\usebox{\tmpfige}\\
notation : &
$0$ & $00$ & $10$ & $11$ & $1$
\end{tabular}
\end{center}
Here, the top $k$-blocks can be anything that comes
between the supporting $1$-block
and the covering $2$-block (inclusive),
but may not be the covering $1$-block.
Columns that are \emph{related} (Definition~\ref{df:35})
are denoted by the same notation.
Notice that we have taken the signature of the corresponding
perfect crystal element for the notation of each Young wall columns.

The following table lists all left-right pairs for which the
signatures to be written under the path description and
Young wall description are not trivially the same.
The signatures in the table body are what should be written as the
left-$\vphi$ and right-$\veps$ signatures under the two Young wall columns.
\begin{center}
\renewcommand{\arraystretch}{1.4}
\begin{tabular}{c*{3}{p{0.4em}c}}
\begin{texdraw}
\fontsize{7}{7}\selectfont
\textref h:C v:C
\drawdim em
\move(0 0)\lvec(-3 1.5)
\htext(-2.7 0.5){left}
\htext(-0.7 1.4){right}
\move(0.3 1.65)
\end{texdraw}&&
$10$ && $11$ && $1$\\
$0$ && $\cdot$ && $1$ && $\cdot\;/\ 01$\\
$00$ && $0$ && $01$ && $0$\\
$10$ && $\cdot$ && $1$ && $\cdot$
\end{tabular}
\end{center}
The case when the left column is $0$ and the right column is $1$
breaks up into two cases.
Depending on the two top $k$-blocks for the left and right columns,
which may be distinct,
the left-$\vphi$ and the right-$\veps$ signatures to be placed
under the two columns could be either nothing or $01$.
In the latter case the signature agrees trivially with that of
the path description.

We can easily see that the signatures agree with that of the
corresponding path description in all the above cases
after $(0,1)$-pair cancellations.
For all other possible left-right pairs not covered in this table,
the left-$\vphi$ and right-$\veps$ signatures to be placed under
the Young wall columns agree exactly with the corresponding path
signatures, that is, even before the $(0,1)$-pair cancellations.

For $0 < i < n$, the following notation will be used.
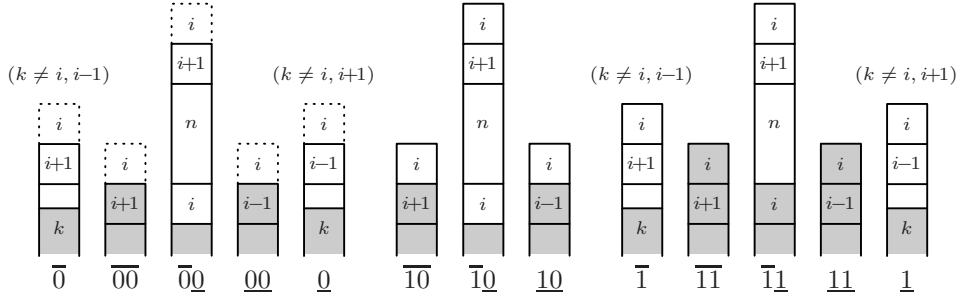
\savebox{\tmpfiga}{\begin{texdraw}
\fontsize{7}{7}\selectfont
\textref h:C v:C
\drawdim em
\setunitscale 1.9
\move(0 -0.75)\lvec(0 1)\lvec(1 1)\lvec(1 -0.75)\ifill f:0.8
\move(0 -0.8)\lvec(0 1)\lvec(1 1)\lvec(1 -0.8)
\move(0 0)\lvec(1 0)
\htext(0.5 0.5){$i\!\!+\!\!1$}
\htext(0.5 1.5){$i$}
\lpatt(0.03 0.17)
\move(0 1)\lvec(0 2)\lvec(1 2)\lvec(1 1)
\end{texdraw}}
\savebox{\tmpfigb}{\begin{texdraw}
\fontsize{7}{7}\selectfont
\textref h:C v:C
\drawdim em
\setunitscale 1.9
\move(0 -0.75)\lvec(0 1)\lvec(1 1)\lvec(1 -0.75)\ifill f:0.8
\move(0 -0.8)\lvec(0 1)\lvec(1 1)\lvec(1 -0.8)
\move(0 0)\lvec(1 0)
\htext(0.5 0.5){$i\!\!-\!\!1$}
\htext(0.5 1.5){$i$}
\lpatt(0.03 0.17)
\move(0 1)\lvec(0 2)\lvec(1 2)\lvec(1 1)
\end{texdraw}}
\savebox{\tmpfigc}{\begin{texdraw}
\fontsize{7}{7}\selectfont
\textref h:C v:C
\drawdim em
\setunitscale 1.9
\move(0 -0.75)\lvec(0 0)\lvec(1 0)\lvec(1 -0.75)\ifill f:0.8
\move(0 -0.8)\lvec(0 4.5)\lvec(1 4.5)\lvec(1 -0.8)
\move(0 0)\lvec(1 0)
\move(0 1)\lvec(1 1)
\move(0 3.5)\lvec(1 3.5)
\htext(0.5 0.5){$i$}
\htext(0.5 4){$i\!\!+\!\!1$}
\htext(0.5 2.5){$n$}
\move(0 4.5)
\bsegment
\lpatt(0.03 0.17)
\move(0 0)\lvec(0 1)\lvec(1 1)\lvec(1 0)
\htext(0.5 0.5){$i$}
\esegment
\end{texdraw}}
\savebox{\tmpfige}{\begin{texdraw}
\fontsize{7}{7}\selectfont
\textref h:C v:C
\drawdim em
\setunitscale 1.9
\move(0 -0.75)\lvec(0 0)\lvec(1 0)\lvec(1 -0.75)\ifill f:0.8
\move(0 -0.8)\lvec(0 4.5)\lvec(1 4.5)\lvec(1 -0.8)
\move(0 0)\lvec(1 0)
\move(0 1)\lvec(1 1)
\move(0 3.5)\lvec(1 3.5)
\htext(0.5 0.5){$i$}
\htext(0.5 4){$i\!\!+\!\!1$}
\htext(0.5 2.5){$n$}
\move(0 4.5)
\bsegment
\move(0 0)\lvec(0 1)\lvec(1 1)\lvec(1 0)
\htext(0.5 0.5){$i$}
\esegment
\end{texdraw}}
\savebox{\tmpfigf}{\begin{texdraw}
\fontsize{7}{7}\selectfont
\textref h:C v:C
\drawdim em
\setunitscale 1.9
\move(0 -0.75)\lvec(0 1)\lvec(1 1)\lvec(1 -0.75)\ifill f:0.8
\move(0 -0.8)\lvec(0 1)\lvec(1 1)\lvec(1 -0.8)
\move(0 0)\lvec(1 0)
\htext(0.5 0.5){$i\!\!+\!\!1$}
\htext(0.5 1.5){$i$}
\move(0 1)\lvec(0 2)\lvec(1 2)\lvec(1 1)
\end{texdraw}}
\savebox{\tmpfigg}{\begin{texdraw}
\fontsize{7}{7}\selectfont
\textref h:C v:C
\drawdim em
\setunitscale 1.9
\move(0 -0.75)\lvec(0 1)\lvec(1 1)\lvec(1 -0.75)\ifill f:0.8
\move(0 -0.8)\lvec(0 1)\lvec(1 1)\lvec(1 -0.8)
\move(0 0)\lvec(1 0)
\htext(0.5 0.5){$i\!\!-\!\!1$}
\htext(0.5 1.5){$i$}
\move(0 1)\lvec(0 2)\lvec(1 2)\lvec(1 1)
\end{texdraw}}
\savebox{\tmpfigh}{\begin{texdraw}
\fontsize{7}{7}\selectfont
\textref h:C v:C
\drawdim em
\setunitscale 1.9
\move(0 -1.75)\lvec(0 -0.6)\lvec(1 -0.6)\lvec(1 -1.75)\ifill f:0.8
\move(0 -0.6)\lvec(1 -0.6)
\htext(0.5 -1.1){$k$}
\move(0 -1.8)\lvec(0 1)\lvec(1 1)\lvec(1 -1.8)
\move(0 0)\lvec(1 0)
\htext(0.5 0.5){$i\!\!+\!\!1$}
\htext(0.5 1.5){$i$}
\lpatt(0.03 0.17)
\move(0 1)\lvec(0 2)\lvec(1 2)\lvec(1 1)
\htext(0.5 2.7){$(k\neq i,i\!\!-\!\!1)$}
\end{texdraw}}
\savebox{\tmpfigi}{\begin{texdraw}
\fontsize{7}{7}\selectfont
\textref h:C v:C
\drawdim em
\setunitscale 1.9
\move(0 -1.75)\lvec(0 -0.6)\lvec(1 -0.6)\lvec(1 -1.75)\ifill f:0.8
\move(0 -0.6)\lvec(1 -0.6)
\htext(0.5 -1.1){$k$}
\move(0 -1.8)\lvec(0 1)\lvec(1 1)\lvec(1 -1.8)
\move(0 0)\lvec(1 0)
\htext(0.5 0.5){$i\!\!-\!\!1$}
\htext(0.5 1.5){$i$}
\lpatt(0.03 0.17)
\move(0 1)\lvec(0 2)\lvec(1 2)\lvec(1 1)
\htext(0.5 2.7){$(k\neq i,i\!\!+\!\!1)$}
\end{texdraw}}
\savebox{\tmpfigj}{\begin{texdraw}
\fontsize{7}{7}\selectfont
\textref h:C v:C
\drawdim em
\setunitscale 1.9
\move(0 -0.75)\lvec(0 2)\lvec(1 2)\lvec(1 -0.75)\ifill f:0.8
\move(0 -0.8)\lvec(0 1)\lvec(1 1)\lvec(1 -0.8)
\move(0 0)\lvec(1 0)
\htext(0.5 0.5){$i\!\!+\!\!1$}
\htext(0.5 1.5){$i$}
\move(0 1)\lvec(0 2)\lvec(1 2)\lvec(1 1)
\end{texdraw}}
\savebox{\tmpfigk}{\begin{texdraw}
\fontsize{7}{7}\selectfont
\textref h:C v:C
\drawdim em
\setunitscale 1.9
\move(0 -0.75)\lvec(0 2)\lvec(1 2)\lvec(1 -0.75)\ifill f:0.8
\move(0 -0.8)\lvec(0 1)\lvec(1 1)\lvec(1 -0.8)
\move(0 0)\lvec(1 0)
\htext(0.5 0.5){$i\!\!-\!\!1$}
\htext(0.5 1.5){$i$}
\move(0 1)\lvec(0 2)\lvec(1 2)\lvec(1 1)
\end{texdraw}}
\savebox{\tmpfigl}{\begin{texdraw}
\fontsize{7}{7}\selectfont
\textref h:C v:C
\drawdim em
\setunitscale 1.9
\move(0 -0.75)\lvec(0 1)\lvec(1 1)\lvec(1 -0.75)\ifill f:0.8
\move(0 -0.8)\lvec(0 4.5)\lvec(1 4.5)\lvec(1 -0.8)
\move(0 0)\lvec(1 0)
\move(0 1)\lvec(1 1)
\move(0 3.5)\lvec(1 3.5)
\htext(0.5 0.5){$i$}
\htext(0.5 4){$i\!\!+\!\!1$}
\htext(0.5 2.5){$n$}
\move(0 4.5)
\bsegment
\move(0 0)\lvec(0 1)\lvec(1 1)\lvec(1 0)
\htext(0.5 0.5){$i$}
\esegment
\end{texdraw}}
\savebox{\tmpfigm}{\begin{texdraw}
\fontsize{7}{7}\selectfont
\textref h:C v:C
\drawdim em
\setunitscale 1.9
\move(0 -1.75)\lvec(0 -0.6)\lvec(1 -0.6)\lvec(1 -1.75)\ifill f:0.8
\move(0 -0.6)\lvec(1 -0.6)
\htext(0.5 -1.1){$k$}
\move(0 -1.8)\lvec(0 1)\lvec(1 1)\lvec(1 -1.8)
\move(0 0)\lvec(1 0)
\htext(0.5 0.5){$i\!\!+\!\!1$}
\htext(0.5 1.5){$i$}
\move(0 1)\lvec(0 2)\lvec(1 2)\lvec(1 1)
\htext(0.5 2.7){$(k\neq i,i\!\!-\!\!1)$}
\end{texdraw}}
\savebox{\tmpfign}{\begin{texdraw}
\fontsize{7}{7}\selectfont
\textref h:C v:C
\drawdim em
\setunitscale 1.9
\move(0 -1.75)\lvec(0 -0.6)\lvec(1 -0.6)\lvec(1 -1.75)\ifill f:0.8
\move(0 -0.6)\lvec(1 -0.6)
\htext(0.5 -1.1){$k$}
\move(0 -1.8)\lvec(0 1)\lvec(1 1)\lvec(1 -1.8)
\move(0 0)\lvec(1 0)
\htext(0.5 0.5){$i\!\!-\!\!1$}
\htext(0.5 1.5){$i$}
\move(0 1)\lvec(0 2)\lvec(1 2)\lvec(1 1)
\htext(0.5 2.7){$(k\neq i,i\!\!+\!\!1)$}
\end{texdraw}}
\begin{center}
\begin{tabular}{ccccc}
\usebox{\tmpfigh}&
\usebox{\tmpfiga}&
\usebox{\tmpfigc}&
\usebox{\tmpfigb}&
\usebox{\tmpfigi}\\
$\overline{0}$ &
$\overline{00}$ &
$\overline{0}\underline{0}$ &
$\underline{00}$ &
$\underline{0}$
\end{tabular}\quad
\begin{tabular}{ccc}
\usebox{\tmpfigf}&
\usebox{\tmpfige}&
\usebox{\tmpfigg}\\
$\overline{10}$ &
$\overline{1}\underline{0}$ &
$\underline{10}$
\end{tabular}\quad
\begin{tabular}{ccccc}
\usebox{\tmpfigm}&
\usebox{\tmpfigj}&
\usebox{\tmpfigl}&
\usebox{\tmpfigk}&
\usebox{\tmpfign}\\
$\overline{1}$ &
$\overline{11}$ &
$\overline{1}\underline{1}$ &
$\underline{11}$ &
$\underline{1}$
\end{tabular}
\end{center}
Again, we have used notation that reflect the signature of
the corresponding perfect crystal element.
Underlines and overlines show whether the $i$-blocks and $i$-slots of
the slice are in supporting or covering positions.
The following table gives the Young wall signatures for all
the nontrivial pairs.
\begin{center}
\renewcommand{\arraystretch}{1.4}
\begin{tabular}{c*{8}{p{0.1em}c}}
\begin{texdraw}
\fontsize{7}{7}\selectfont
\textref h:C v:C
\drawdim em
\move(0 0)\lvec(-3 1.5)
\htext(-2.7 0.5){left}
\htext(-0.7 1.4){right}
\move(0.3 1.65)
\end{texdraw}&&
$\overline{10}$ &&
$\overline{1}\underline{0}$ &&
$\underline{10}$ &&
$\overline{1}$ &&
$\overline{11}$ &&
$\overline{1}\underline{1}$ &&
$\underline{11}$ &&
$\underline{1}$\\
$\overline{0}$
  && $\cdot$ && $\cdot\;/\,01$ && && $\cdot\;/\,01$ && $1$ && $1\,/\,011$\\
$\overline{00}$
  && $0$ && $0$ && && $0$ && $\cdot$ && $01$\\
$\overline{0}\underline{0}$
  && $0$ && $0$ && $0$ && $0\,/\,001$ && $01$ && $01$ && $01$ && $0\,/\,001$\\
$\underline{00}$
  && && $0$ && $0$ && && && $01$ && $\cdot$ && $0$\\
$\underline{0}$
  && && $\cdot\;/\,01$ && $\cdot$ && && && $1\;/\,011$ && $1$ && $\cdot\;/\,01$\\
$\overline{10}$
  && $\cdot$ && $\cdot$ && && $\cdot$ && $1$ && $1$\\
$\overline{1}\underline{0}$
  && $\cdot$ && $\cdot$ && $\cdot$ && $\cdot\;/\,01$ && $1$ && $1$ && $1$
  && $\cdot\;/\,01$\\
$\underline{10}$
  && && $\cdot$ && $\cdot$ && && && $1$ && $1$ && $\cdot$
\end{tabular}
\end{center}
The blank slots are where the signatures agree trivially.
As before, many of the cases break up into subcases, depending on the
top $k$-blocks for $\underline{0}$, $\overline{0}$, $\underline{1}$,
and $\overline{1}$.

For the remaining $i=n$ case, the following notation will be used.
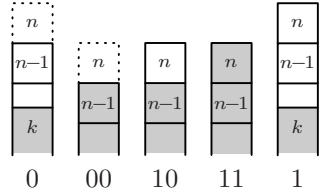
\savebox{\tmpfiga}{\begin{texdraw}
\fontsize{7}{7}\selectfont
\textref h:C v:C
\drawdim em
\setunitscale 1.9
\move(0 -0.75)\lvec(0 1)\lvec(1 1)\lvec(1 -0.75)\ifill f:0.8
\move(0 -0.8)\lvec(0 1)\lvec(1 1)\lvec(1 -0.8)
\move(0 0)\lvec(1 0)
\htext(0.5 0.5){$n\!\!-\!\!1$}
\htext(0.5 1.5){$n$}
\lpatt(0.03 0.17)
\move(0 1)\lvec(0 2)\lvec(1 2)\lvec(1 1)
\end{texdraw}}
\savebox{\tmpfigf}{\begin{texdraw}
\fontsize{7}{7}\selectfont
\textref h:C v:C
\drawdim em
\setunitscale 1.9
\move(0 -0.75)\lvec(0 1)\lvec(1 1)\lvec(1 -0.75)\ifill f:0.8
\move(0 -0.8)\lvec(0 1)\lvec(1 1)\lvec(1 -0.8)
\move(0 0)\lvec(1 0)
\htext(0.5 0.5){$n\!\!-\!\!1$}
\htext(0.5 1.5){$n$}
\move(0 1)\lvec(0 2)\lvec(1 2)\lvec(1 1)
\end{texdraw}}
\savebox{\tmpfigj}{\begin{texdraw}
\fontsize{7}{7}\selectfont
\textref h:C v:C
\drawdim em
\setunitscale 1.9
\move(0 -0.75)\lvec(0 2)\lvec(1 2)\lvec(1 -0.75)\ifill f:0.8
\move(0 -0.8)\lvec(0 1)\lvec(1 1)\lvec(1 -0.8)
\move(0 0)\lvec(1 0)
\htext(0.5 0.5){$n\!\!-\!\!1$}
\htext(0.5 1.5){$n$}
\move(0 1)\lvec(0 2)\lvec(1 2)\lvec(1 1)
\end{texdraw}}
\savebox{\tmpfigi}{\begin{texdraw}
\fontsize{7}{7}\selectfont
\textref h:C v:C
\drawdim em
\setunitscale 1.9
\move(0 -1.75)\lvec(0 -0.6)\lvec(1 -0.6)\lvec(1 -1.75)\ifill f:0.8
\move(0 -0.6)\lvec(1 -0.6)
\htext(0.5 -1.1){$k$}
\move(0 -1.8)\lvec(0 1)\lvec(1 1)\lvec(1 -1.8)
\move(0 0)\lvec(1 0)
\htext(0.5 0.5){$n\!\!-\!\!1$}
\htext(0.5 1.5){$n$}
\lpatt(0.03 0.17)
\move(0 1)\lvec(0 2)\lvec(1 2)\lvec(1 1)
\end{texdraw}}
\savebox{\tmpfign}{\begin{texdraw}
\fontsize{7}{7}\selectfont
\textref h:C v:C
\drawdim em
\setunitscale 1.9
\move(0 -1.75)\lvec(0 -0.6)\lvec(1 -0.6)\lvec(1 -1.75)\ifill f:0.8
\move(0 -0.6)\lvec(1 -0.6)
\htext(0.5 -1.1){$k$}
\move(0 -1.8)\lvec(0 1)\lvec(1 1)\lvec(1 -1.8)
\move(0 0)\lvec(1 0)
\htext(0.5 0.5){$n\!\!-\!\!1$}
\htext(0.5 1.5){$n$}
\move(0 1)\lvec(0 2)\lvec(1 2)\lvec(1 1)
\end{texdraw}}
\begin{center}
\begin{tabular}{ccccc}
\usebox{\tmpfigi}&
\usebox{\tmpfiga}&
\usebox{\tmpfigf}&
\usebox{\tmpfigj}&
\usebox{\tmpfign}\\
$0$ & $00$ & $10$ & $11$ & $1$
\end{tabular}
\end{center}
Here, the top $k$-blocks for $0$ and $1$
may be taken to be anything except for the $n$-block and
the supporting $(n-1)$-block.
It could also be two $0$-blocks.
The following table lists all nontrivial pairs.
\begin{center}
\renewcommand{\arraystretch}{1.4}
\begin{tabular}{c*{3}{p{0.4em}c}}
\begin{texdraw}
\fontsize{7}{7}\selectfont
\textref h:C v:C
\drawdim em
\move(0 0)\lvec(-3 1.5)
\htext(-2.7 0.5){left}
\htext(-0.7 1.4){right}
\move(0.3 1.65)
\end{texdraw}&&
$10$ && $11$ && $1$\\
$0$ && $\cdot$ && $1$ && $\cdot\;/\ 01$\\
$00$ && $0$ && $01$ && $0$\\
$10$ && $\cdot$ && $1$ && $\cdot$
\end{tabular}
\end{center}

This completes the proof of Lemma~\ref{lem:55} and of Theorem~\ref{thm:54}.

We close this paper with an example.
\begin{example}
Following is the top part
of the crystal graph $\rpwspace(\La_0)$ for $\uq(C_2^{(1)})$.
\begin{center}
\begin{texdraw}
\drawdim in \arrowheadsize l:0.065 w:0.03 \arrowheadtype t:F
\fontsize{5}{5}\selectfont
\textref h:C v:C
\drawdim em
\setunitscale 1.9
\move(-13.6 0.32)
\bsegment
\setsegscale 0.7
\move(-1.25 0)\lvec(0 0)\lvec(0 0.5)\lvec(-1.25 0.5)\lvec(-1.25 0)\ifill f:0.8
\move(-1.3 0)\lvec(0 0)\lvec(0 0.5)\lvec(-1.3 0.5)
\move(-1 0)\lvec(-1 0.5)
\htext(-0.5 0.27){$0$}
\esegment
\move(-9.9 0.25)
\bsegment
\setsegscale 0.7
\move(-1.25 0)\lvec(0 0)\lvec(0 1)\lvec(-1 1)\lvec(-1 0.5)
\lvec(-1.25 0.5)\lvec(-1.25 0)\ifill f:0.8
\move(-1.3 0)\lvec(0 0)\lvec(0 1)\lvec(-1 1)\lvec(-1 0)
\move(-1.3 0.5)\lvec(0 0.5)
\htext(-0.5 0.27){$0$}
\htext(-0.5 0.77){$0$}
\esegment
\move(-6.2 0)
\bsegment
\setsegscale 0.7
\move(-1.25 0)\lvec(0 0)\lvec(0 1)\lvec(-1 1)\lvec(-1 0.5)
\lvec(-1.25 0.5)\lvec(-1.25 0)\ifill f:0.8
\move(-1.3 0)\lvec(0 0)\lvec(0 2)\lvec(-1 2)\lvec(-1 0)
\move(-1.3 0.5)\lvec(0 0.5)
\move(-1 1)\lvec(0 1)
\htext(-0.5 0.27){$0$}
\htext(-0.5 0.77){$0$}
\htext(-0.5 1.5){$1$}
\esegment
\move(-2.5 0)
\bsegment
\move(0 -1.5)
\bsegment
\setsegscale 0.7
\move(-1.25 0)\lvec(0 0)\lvec(0 2)\lvec(-1 2)\lvec(-1 0.5)
\lvec(-1.25 0.5)\lvec(-1.25 0)\ifill f:0.8
\move(-1.3 0)\lvec(0 0)\lvec(0 2)\lvec(-1 2)\lvec(-1 0)
\move(-1.3 0.5)\lvec(0 0.5)
\move(-1 1)\lvec(0 1)
\htext(-0.5 0.27){$0$}
\htext(-0.5 0.77){$0$}
\htext(-0.5 1.5){$1$}
\esegment
\move(0 1.5)
\bsegment
\setsegscale 0.7
\move(-1.25 0)\lvec(0 0)\lvec(0 1)\lvec(-1 1)\lvec(-1 0.5)
\lvec(-1.25 0.5)\lvec(-1.25 0)\ifill f:0.8
\move(-1.3 0)\lvec(0 0)\lvec(0 3)\lvec(-1 3)\lvec(-1 0)
\move(-1.3 0.5)\lvec(0 0.5)
\move(-1 1)\lvec(0 1)
\move(-1 2)\lvec(0 2)
\htext(-0.5 0.27){$0$}
\htext(-0.5 0.77){$0$}
\htext(-0.5 1.5){$1$}
\htext(-0.5 2.5){$2$}
\esegment
\esegment
\move(2 0)
\bsegment
\move(0 -3)
\bsegment
\setsegscale 0.7
\move(-2.25 0)\lvec(0 0)\lvec(0 2)\lvec(-1 2)\lvec(-1 1)\lvec(-2 1)
\lvec(-2 0.5)\lvec(-2.25 0.5)\lvec(-2.25 0)\ifill f:0.8
\move(-2.3 0)\lvec(0 0)\lvec(0 2)\lvec(-1 2)\lvec(-1 0)
\move(-2.3 0.5)\lvec(0 0.5)
\move(-2 0)\lvec(-2 1)\lvec(0 1)
\htext(-1.5 0.27){$0$}
\htext(-1.5 0.77){$0$}
\htext(-0.5 0.27){$0$}
\htext(-0.5 0.77){$0$}
\htext(-0.5 1.5){$1$}
\esegment
\move(0 -0.3)
\bsegment
\setsegscale 0.7
\move(-2.25 0)\lvec(0 0)\lvec(0 2)\lvec(-1 2)\lvec(-1 0.5)
\lvec(-2.25 0.5)\lvec(-2.25 0)\ifill f:0.8
\move(-2.3 0)\lvec(0 0)\lvec(0 3)\lvec(-1 3)\lvec(-1 0)
\move(-2.3 0.5)\lvec(0 0.5)
\move(-2 0)\lvec(-2 0.5)
\move(-1 1)\lvec(0 1)
\move(-1 2)\lvec(0 2)
\htext(-1.5 0.27){$0$}
\htext(-0.5 0.27){$0$}
\htext(-0.5 0.77){$0$}
\htext(-0.5 1.5){$1$}
\htext(-0.5 2.5){$2$}
\esegment
\move(0 3)
\bsegment
\setsegscale 0.7
\move(-2.25 0)\lvec(0 0)\lvec(0 1)\lvec(-1 1)\lvec(-1 0.5)
\lvec(-2.25 0.5)\lvec(-2.25 0)\ifill f:0.8
\move(-2.3 0)\lvec(0 0)\lvec(0 4)\lvec(-1 4)\lvec(-1 0)
\move(-2.3 0.5)\lvec(0 0.5)
\move(-2 0)\lvec(-2 0.5)
\move(-1 1)\lvec(0 1)
\move(-1 2)\lvec(0 2)
\move(-1 3)\lvec(0 3)
\htext(-1.5 0.27){$0$}
\htext(-0.5 0.27){$0$}
\htext(-0.5 0.77){$0$}
\htext(-0.5 1.5){$1$}
\htext(-0.5 2.5){$2$}
\htext(-0.5 3.5){$1$}
\esegment
\esegment
\move(6.5 -0.5)
\bsegment
\move(0 -4.5)
\bsegment
\setsegscale 0.7
\move(-2.25 0)\lvec(0 0)\lvec(0 2)\lvec(-1 2)\lvec(-1 1)\lvec(-2 1)
\lvec(-2 0.5)\lvec(-2.25 0.5)\lvec(-2.25 0)\ifill f:0.8
\move(-2.3 0)\lvec(0 0)\lvec(0 3)\lvec(-1 3)\lvec(-1 0)
\move(-2.3 0.5)\lvec(0 0.5)
\move(-2 0)\lvec(-2 1)\lvec(0 1)
\move(-1 2)\lvec(0 2)
\htext(-1.5 0.27){$0$}
\htext(-1.5 0.77){$0$}
\htext(-0.5 0.27){$0$}
\htext(-0.5 0.77){$0$}
\htext(-0.5 1.5){$1$}
\htext(-0.5 2.5){$2$}
\esegment
\move(0 -1.35)
\bsegment
\setsegscale 0.7
\move(-2.25 0)\lvec(0 0)\lvec(0 3)\lvec(-1 3)\lvec(-1 0.5)
\lvec(-2.25 0.5)\lvec(-2.25 0)\ifill f:0.8
\move(-2.3 0)\lvec(0 0)\lvec(0 3)\lvec(-1 3)\lvec(-1 0)
\move(-2.3 0.5)\lvec(0 0.5)
\move(-2 0)\lvec(-2 0.5)
\move(-1 1)\lvec(0 1)
\move(-1 2)\lvec(0 2)
\htext(-1.5 0.27){$0$}
\htext(-0.5 0.27){$0$}
\htext(-0.5 0.77){$0$}
\htext(-0.5 1.5){$1$}
\htext(-0.5 2.5){$2$}
\esegment
\move(0 1.8)
\bsegment
\setsegscale 0.7
\move(-2.25 0)\lvec(0 0)\lvec(0 1)\lvec(-2 1)
\lvec(-2 0.5)\lvec(-2.25 0.5)\lvec(-2.25 0)\ifill f:0.8
\move(-2.3 0)\lvec(0 0)\lvec(0 4)\lvec(-1 4)\lvec(-1 0)
\move(-2.3 0.5)\lvec(0 0.5)
\move(-2 0)\lvec(-2 1)\lvec(0 1)
\move(-1 2)\lvec(0 2)
\move(-1 3)\lvec(0 3)
\htext(-1.5 0.27){$0$}
\htext(-1.5 0.77){$0$}
\htext(-0.5 0.27){$0$}
\htext(-0.5 0.77){$0$}
\htext(-0.5 1.5){$1$}
\htext(-0.5 2.5){$2$}
\htext(-0.5 3.5){$1$}
\esegment
\move(0 5.4)
\bsegment
\setsegscale 0.7
\move(-2.25 0)\lvec(0 0)\lvec(0 2)\lvec(-1 2)\lvec(-1 0.5)
\lvec(-2.25 0.5)\lvec(-2.25 0)\ifill f:0.8
\move(-2.3 0)\lvec(0 0)\lvec(0 4)\lvec(-1 4)\lvec(-1 0)
\move(-2.3 0.5)\lvec(0 0.5)
\move(-2 0)\lvec(-2 0.5)
\move(-1 1)\lvec(0 1)
\move(-1 2)\lvec(0 2)
\move(-1 3)\lvec(0 3)
\htext(-1.5 0.27){$0$}
\htext(-0.5 0.27){$0$}
\htext(-0.5 0.77){$0$}
\htext(-0.5 1.5){$1$}
\htext(-0.5 2.5){$2$}
\htext(-0.5 3.5){$1$}
\esegment
\esegment
\move(11 -1)
\bsegment
\move(0 -6)
\bsegment
\setsegscale 0.7
\move(-2.25 0)\lvec(0 0)\lvec(0 2)\lvec(-1 2)\lvec(-1 1)\lvec(-2 1)
\lvec(-2 0.5)\lvec(-2.25 0.5)\lvec(-2.25 0)\ifill f:0.8
\move(-2.3 0)\lvec(0 0)\lvec(0 3)\lvec(-1 3)\lvec(-1 0)
\move(-2 0)\lvec(-2 2)\lvec(0 2)
\move(-2.3 0.5)\lvec(0 0.5)
\move(-2 1)\lvec(0 1)
\htext(-1.5 0.27){$0$}
\htext(-1.5 0.77){$0$}
\htext(-0.5 0.27){$0$}
\htext(-0.5 0.77){$0$}
\htext(-1.5 1.5){$1$}
\htext(-0.5 1.5){$1$}
\htext(-0.5 2.5){$2$}
\esegment
\move(0 -3)
\bsegment
\setsegscale 0.7
\move(-2.25 0)\lvec(0 0)\lvec(0 3)\lvec(-1 3)\lvec(-1 1)\lvec(-2 1)
\lvec(-2 0.5)\lvec(-2.25 0.5)\lvec(-2.25 0)\ifill f:0.8
\move(-2.3 0)\lvec(0 0)\lvec(0 3)\lvec(-1 3)\lvec(-1 0)
\move(-2.3 0.5)\lvec(0 0.5)
\move(-2 0)\lvec(-2 1)\lvec(0 1)
\move(-1 2)\lvec(0 2)
\htext(-1.5 0.27){$0$}
\htext(-1.5 0.77){$0$}
\htext(-0.5 0.27){$0$}
\htext(-0.5 0.77){$0$}
\htext(-0.5 1.5){$1$}
\htext(-0.5 2.5){$2$}
\esegment
\move(0 0)
\bsegment
\setsegscale 0.7
\move(-2.25 0)\lvec(0 0)\lvec(0 3)\lvec(-1 3)\lvec(-1 0.5)
\lvec(-2.25 0.5)\lvec(-2.25 0)\ifill f:0.8
\move(-2.3 0)\lvec(0 0)\lvec(0 4)\lvec(-1 4)\lvec(-1 0)
\move(-2.3 0.5)\lvec(0 0.5)
\move(-2 0)\lvec(-2 0.5)
\move(-1 1)\lvec(0 1)
\move(-1 2)\lvec(0 2)
\move(-1 3)\lvec(0 3)
\htext(-1.5 0.27){$0$}
\htext(-0.5 0.27){$0$}
\htext(-0.5 0.77){$0$}
\htext(-0.5 1.5){$1$}
\htext(-0.5 2.5){$2$}
\htext(-0.5 3.5){$1$}
\esegment
\move(0 3.7)
\bsegment
\setsegscale 0.7
\move(-2.25 0)\lvec(0 0)\lvec(0 1)\lvec(-2 1)
\lvec(-2 0.5)\lvec(-2.25 0.5)\lvec(-2.25 0)\ifill f:0.8
\move(-2.3 0)\lvec(0 0)\lvec(0 4)\lvec(-1 4)\lvec(-1 0)
\move(-2 0)\lvec(-2 2)\lvec(0 2)
\move(-2.3 0.5)\lvec(0 0.5)
\move(-2 1)\lvec(0 1)
\move(-1 3)\lvec(0 3)
\htext(-1.5 0.27){$0$}
\htext(-1.5 0.77){$0$}
\htext(-0.5 0.27){$0$}
\htext(-0.5 0.77){$0$}
\htext(-1.5 1.5){$1$}
\htext(-0.5 1.5){$1$}
\htext(-0.5 2.5){$2$}
\htext(-0.5 3.5){$1$}
\esegment
\move(0 7.4)
\bsegment
\setsegscale 0.7
\move(-2.25 0)\lvec(0 0)\lvec(0 2)\lvec(-1 2)\lvec(-1 1)\lvec(-2 1)
\lvec(-2 0.5)\lvec(-2.25 0.5)\lvec(-2.25 0)\ifill f:0.8
\move(-2.3 0)\lvec(0 0)\lvec(0 4)\lvec(-1 4)\lvec(-1 0)
\move(-2.3 0.5)\lvec(0 0.5)
\move(-2 0)\lvec(-2 1)\lvec(0 1)
\move(-1 2)\lvec(0 2)
\move(-1 3)\lvec(0 3)
\htext(-1.5 0.27){$0$}
\htext(-1.5 0.77){$0$}
\htext(-0.5 0.27){$0$}
\htext(-0.5 0.77){$0$}
\htext(-0.5 1.5){$1$}
\htext(-0.5 2.5){$2$}
\htext(-0.5 3.5){$1$}
\esegment
\esegment
\move(-13.2 0.5)\ravec(2 0)\htext(-12.3 0.8){$0$}
\move(-9.5 0.5)\ravec(2 0)\htext(-8.6 0.8){$1$}
\move(-5.8 0.9)\ravec(2 1)\htext(-4.9 1.75){$2$}
\move(-5.8 0.4)\ravec(2 -1)\htext(-4.9 0.3){$1$}
\move(-2.1 2.7)\ravec(2 1)\htext(-1.2 3.6){$1$}
\move(-2.1 -0.6)\ravec(2 1)\htext(-1.2 0.25){$2$}
\move(-2.1 -1.1)\ravec(2 -1)\htext(-1.2 -1.2){$0$}
\move(2.5 4.6)\ravec(2 1)\htext(3.4 5.9){$1$}
\move(2.5 3.8)\ravec(2 -1)\htext(3.4 3.7){$0$}
\move(2.5 0.6)\ravec(2 -1)\htext(3.4 0.6){$2$}
\move(2.5 -0.3)\ravec(2 -2.5)\htext(3.4 -0.85){$0$}
\move(2.5 -2.9)\ravec(2 -1)\htext(3.4 -2.95){$2$}
\move(6.9 6.4)\ravec(2 1)\htext(7.8 7.3){$0$}
\move(6.9 5.3)\ravec(2 -3.8)\htext(7.8 4.3){$2$}
\move(6.9 2.5)\ravec(2 1)\htext(7.2 3){$1$}
\move(6.9 -0.9)\ravec(2 1)\htext(7.8 -0.05){$1$}
\move(6.9 -1.5)\ravec(2 -1)\htext(7.8 -1.55){$0$}
\move(6.9 -4.3)\ravec(2 1)\htext(7.8 -3.45){$2$}
\move(6.9 -4.9)\ravec(2 -1)\htext(7.8 -5){$1$}
\htext(12 7.4){$\cdots$}
\htext(12 3.7){$\cdots$}
\htext(12 0){$\cdots$}
\htext(12 -3.45){$\cdots$}
\htext(12 -6.5){$\cdots$}
\end{texdraw}
\end{center}
\end{example}

\bibliographystyle{amsplain}

\def\cprime{$'$}
\providecommand{\bysame}{\leavevmode\hbox to3em{\hrulefill}\thinspace}
\providecommand{\MR}{\relax\ifhmode\unskip\space\fi MR }
\providecommand{\MRhref}[2]{%
  \href{http://www.ams.org/mathscinet-getitem?mr=#1}{#2}
}
\providecommand{\href}[2]{#2}

\end{document}